%% file: main.tex
\documentclass[12pt]{amsart}

\input{preamble.tex}

\title[Goal-Oriented A-Posteriori Estimation of Model Error]{Goal-Oriented A-Posteriori Estimation of Model Error as an Aid to Parameter Estimation}

\author{Prashant K. Jha and J. Tinsley Oden}
\address{Oden Institute for Computational Engineering and Sciences, The University of Texas at Austin, Austin, TX 78712}
\email{prashant.jha@austin.utexas.edu, oden@oden.utexas.edu}

\begin{document}

\begin{abstract}
In this work, a Bayesian model calibration framework is presented that utilizes goal-oriented a-posterior error estimates in quantities of interest (QoIs) for classes of high-fidelity models characterized by PDEs. It is shown that for a large class of computational models, it is possible to develop a computationally inexpensive procedure for calibrating parameters of high-fidelity models of physical events when the parameters of low-fidelity (surrogate) models are known with acceptable accuracy. The main ingredients in the proposed model calibration scheme are goal-oriented a-posteriori estimates of error in QoIs computed using a so-called lower fidelity model compared to those of an uncalibrated higher fidelity model. The estimates of error in QoIs are used to define likelihood functions in Bayesian inversion analysis. A standard Bayesian approach is employed to compute the posterior distribution of model parameters of high-fidelity models. As applications, parameters in a quasi-linear second-order elliptic boundary-value problem (BVP) are calibrated using a second-order linear elliptic BVP. In a second application, parameters of a tumor growth model involving nonlinear time-dependent PDEs are calibrated using a lower fidelity linear tumor growth model with known parameter values.
\end{abstract}

\maketitle

\section{Introduction}
In this work, a fundamental question in predictive computational science embodied in the following scenario is addressed:
\begin{itemize}
\item[1.] Given a high-fidelity model of a class of physical phenomena which possibly involves a large number of unknown or poorly specified parameters; and
\item[2.] Given a lower-fidelity model (or a class of ``surrogate" models) of the same physical events which involves fewer parameters but for which the parameters are known with acceptable precision; and
\item[3.] Further, suppose it is possible to derive \textit{a posterior error estimates} in key \textit{quantities-of-interest}, the QoIs, (so-called ``goal-oriented" estimates) estimating the error in the predictions of the low-fidelity model compared to those of the high-fidelity model.   
\end{itemize}
\textit{Then, can one use such estimates and calibrated lower fidelity models to infer values of the parameters of the high-fidelity model?} Further, can such parameter estimates be made in the presence of uncertainties?

The issue of developing a-posteriori estimates of modeling error in QoIs was taken up in \cite{oden2001goal, oden2002estimation, van2011goal, prudhomme1999goal, prudhomme2003computable, rannacher1997feed}, and generally reduces to estimates of the form,
\begin{equation}\label{eq:errQoI}
	\calQ(u) - \calQ(u_0) = \calR(u_0; p) + r \approx \calR(u_0; p),
\end{equation}
where $\calQ(\cdot)$ is the value of the QoI-functional on solutions $u$ of the forward problem, $u_0$ is an approximation of $u$ provided as a solution of a computable surrogate or reduced-fidelity model, $\calR(u_0; \cdot)$ is one of several possible forms of a residual functional defining the misfit of the surrogate solution to the forward problem, $p$ is the solution of the adjoint problem associated with the forward problem operators and the QoI functional. If $p_0$ is an approximation of $p$ possibly obtained by solving the dual problem associated with the surrogate and the QoI functional then $r = r(u, p, u_0, p_0)$ is the remainder term, often assumed to be negligible, as it involves second or higher derivatives of semilinear forms $\calB$ defining the forward problem and the QoI functional $\calQ$. We remark that several different forms of the residual $\calR(\cdot; \cdot)$ can be considered depending on various simplifying assumptions. 

In this exposition, the use of such error estimates as a source of data in Bayesian framework for parameter estimation in the presence of uncertainties is explored. However, the success of such inference methods, not surprisingly, depends upon the magnitude of the errors, $e_0 = u - u_0$, in the surrogate approximation of the forward problem and the error $\eps_0 = p - p_0$ in the surrogate approximation of the adjoint problem. For sufficiently small errors $(e_0, \eps_0)$, approximate residuals can be derived which employ approximations of these error components determined by solutions of an auxiliary pair of linear variational problems. In such cases, it is argued that quite accurate estimates of key parameters of the high-fidelity model can be obtained at reasonable computational costs. Applications to a representative of a class of quasi-linear elliptic boundary-value problems and to a class of time-dependent models of tumor growth at the tissue scale are presented. Details on the implementation and behavior of the Bayesian methodology for parameter estimation are given for these example applications.

Several studies have been published in recent literature that make use of reduced order models in a Bayesian context designed to accelerate the parameter estimation, e.g., \cite{yan2019adaptive, manzoni2016accurate, yan2019adaptive2pub, conrad2016accelerating, li2018model, biros2011large, lassila2013reduced, galbally2010non, marzouk2009dimensionality, roderick2014proper, frangos2010surrogate}. The goal of this work, however, is different. This work proposes techniques that extract data for Bayesian parameter inversion, not from experimental validation scenarios but by drawing data from low-fidelity surrogates on specific calibrated QoIs and exploiting rigorous a-posterior estimates of the goal-oriented error to infer parameters for a different higher fidelity model. At the same time, the connection between the proposed work and inference problems in which surrogates are used instead of higher-fidelity models can be seen as follows: goal-oriented a-posterior error estimates can be used to correct errors in Bayesian inference calculations in scenarios in which the higher-fidelity model is replaced with surrogates. 

Following this introduction, a brief review of the theory of goal-oriented a-posteriori estimates of modeling error following \cite{oden2001goal, oden2002estimation} is provided and new methods of approximating the pair $(e_0, \eps_0)$ useful to the evaluation of estimates of parameters of the high-fidelity model are presented. In \cref{s:bayes}, a Bayesian inversion framework for estimating parameters of the high-fidelity model in the presence of uncertainty using the lower-fidelity solution and the computable estimates of the QoI error is presented. Representative applications are discussed in \cref{s:examples} together with details on the performance of Markov Chain Monte Carlo (MCMC) calculations used for parameter estimations. Concluding comments are collected in \cref{s:conclusion}. Codes used to obtain the numerical results are available in this link: \url{https://github.com/prashjha/GoalOrientedModelCalibration}.

%%%%%%%%%%%%%%%%%%%%%%%%%%%%%%%%%%%%%%%%%%%%%%%%%%%%%%%%%%%%%%%%%
%%%%%%%%%%%%%%%%%%%%%%%%%%%%%%%%%%%%%%%%%%%%%%%%%%%%%%%%%%%%%%%%%
\section{Goal-Oriented Estimation of the Modeling Error}\label{s:goal}
Consider the abstract nonlinear problem: Find $u \in \calV$ such that
\begin{equation}\label{eq:fwdProbWeak}
	\calB(u; v) = \calF(v), \qquad \forall v \in \calV,
\end{equation}
where $\calB(\cdot; \cdot)$ is a differentiable semilinear form on a Banach space $\calV$, which is linear in arguments following the semicolon but possibly nonlinear in $u$, and $\calF(\cdot)$ is a given linear functional on $\calV$. The primary goal in formulating and solving \eqref{eq:fwdProbWeak} is to determine features of the solution characterized by another possibly nonlinear functional $\calQ$ defined on $\calV$, the quantity of interest. 
At the outset, it is assumed that \eqref{eq:fwdProbWeak} admits a unique solution $u \in \calV$ and that $\calB(\cdot; \cdot)$ and $\calQ(\cdot)$ are differentiable in the G\^ateaux or functional sense to a high degree, generally three or more; that is, the limits such as
\begin{align}
\calB'(u; v, p) &= \lim_{\eta \to 0} \eta^{-1} \sqbrack{\calB(u+\eta v; p) - \calB(u; p)}, \notag \\
\calB''(u; q, v, p) &= \lim_{\eta \to 0} \eta^{-1} \sqbrack{\calB'(u+\eta q; v, p) - \calB'(u; v, p)}, \notag \\
\calB'''(u; r, q, v, p) &= \lim_{\eta \to 0} \eta^{-1} \sqbrack{\calB''(u+\eta r;  q, v, p) - \calB''(u; q, v, p)}, ...
\end{align}
and
\begin{align}
\calQ'(u; v) &= \lim_{\eta \to 0} \eta^{-1} \sqbrack{\calQ(u+\eta v) - \calQ(u)}, \notag \\
\calQ''(u; q, v) &= \lim_{\eta \to 0} \eta^{-1} \sqbrack{\calQ'(u+\eta q; v) - \calQ'(u; v)}, \notag \\
\calQ'''(u; r, q, v) &= \lim_{\eta \to 0} \eta^{-1} \sqbrack{\calQ''(u+\eta r; q, v) - \calQ''(u; q, v)}, ...
\end{align}
exist, $\eta \in \bbR^+$. Following \cite{becker2001optimal, bangerth2013adaptive, oden2001goal, oden2002estimation, prudhomme1999goal, prudhomme2003computable}, it is useful to note that the value of the QoI, $\calQ(u)$, $u$ being the solution of the forward problem \eqref{eq:fwdProbWeak}, can be computed as the solution of the following constrained minimization problem: Find $u \in \calV$ such that
\begin{equation}\label{eq:optProb}
	\calQ(u) = \inf_{v \in M} \calQ(v),
\end{equation}
where
\begin{equation*}
M = \cubrack{v\in \calV: \calB(v; q) = \calF(q), \, \forall q\in \calV}.
\end{equation*}
The constrained optimization problem in \eqref{eq:optProb} the constraint being that the admissible functions $v$ must satisfy the forward problem, i.e., $v$ is such that $\calB(v; q) = \calF(q)$ for all $q\in \calV$, has an associated Lagrangian $L = L(v, q)$ given by
\begin{equation}
L(v, q) = \calQ(v) + \calF(q) - \calB(v; q).
\end{equation}
The constrained optimization problem \eqref{eq:optProb} can now be solved by finding the extremum (critical points) $(u, p)$ of the Lagrangian $L$. 
%The minimizer $u$ corresponds to the first component of a saddle point $(u,p) \in \calV\times \calV$ of the Lagrangian, 
%\begin{equation}
%L(u, p) = \calQ(u) + \calF(p) - \calB(u; p)
%\end{equation}
%with $p$ the Lagrange multiplier, or adjoint or dual variable corresponding to the choice $\calQ$ of the quantity of interest. 
The critical points $(u,p)$ of $L$ are such that $L'((u,p); (v,q)) = 0$, $\forall (v,q) \in \calV\times \calV$, which are solutions of the equations:
\begin{equation}\label{eq:fwdAdj}
	\begin{split}
		\calB(u; q) &= \calF(q), \qquad \forall q \in \calV, \\
		\calB'(u; v, p) &= \calQ'(u; v), \qquad \forall v \in \calV .
	\end{split}
\end{equation}
The first equation in \eqref{eq:fwdAdj} is recognized as the primal or forward problem \eqref{eq:fwdProbWeak} while the second equation is the adjoint or dual problem for $p$ with $u$ specified. The adjoint problem \eqref{eq:fwdAdj}$_2$ is thus a linear variational (weak) formulation for $p$ given $u$ and the quantity of interest functional $\calQ$. 

Let us now suppose that \eqref{eq:fwdProbWeak} is intractable for practical purposes so that we are led to consider a different semilinear form $\calB_0(\cdot; \cdot)$ on $\calV \times \calV$ that may be a coarser lower-fidelity model of the same physical event modeled by \eqref{eq:fwdProbWeak} with solutions $u_0\in \calV_0 \subseteq \calV$, $\calV_0$ being subspace of $\calV$. Thus, a lower-fidelity problem is given by: Find $u_0 \in \calV_0$ such that
\begin{equation}\label{eq:fwdProbLFWeak}
	\calB_0(u_0; v) = \calF(v), \qquad \forall v \in \calV_0 .
\end{equation}
Following the same steps used to obtain \eqref{eq:fwdAdj}, a constrained optimization problem for the lower-fidelity surrogate model \eqref{eq:fwdProbLFWeak} can be formulated to arrive at these surrogate forward and adjoint pair of equations, 
\begin{equation}\label{eq:fwdAdjLF}
	\begin{split}
		\calB_0(u_0; q) &= \calF(q), \qquad \forall q \in \calV_0 \\
		\calB'_0(u_0; v, p_0) &= \calQ'(u_0; v), \qquad \forall v \in \calV_0 .
	\end{split}
\end{equation}
Thus, again, $\calB_0(\cdot ; \cdot)$ is assumed to be differentiable. In many cases, $\calV_0 = \calV$ as \eqref{eq:fwdAdj} and \eqref{eq:fwdAdjLF} are different models of the same physical events. 

The primal and adjoint errors $(e_0, \eps_0)$ are defined as
\begin{equation}\label{eq:err}
e_0 = u - u_0 \qquad \qquad \text{and} \qquad \qquad \eps_0 = p - p_0 .
\end{equation}
The degree to which the reduced model solutions $(u_0, p_0)$ fail to satisfy \eqref{eq:fwdAdj} is characterized by the residuals, $\calR(\cdot; \cdot)$ and $\bar{\calR}(\cdot; \cdot, \cdot)$, defined as:
\begin{equation}\label{eq:residuals}
	\begin{split}
		\calR(u_0; q) &= \calF(q) - \calB(u_0; q) , \qquad \forall q \in \calV, \\
		\bar{\calR}(u_0; v, p_0) &= \calQ'(u_0; v) - \calB'(u_0; v, p_0), \qquad \forall v \in \calV .
	\end{split}
\end{equation}
While the fine and coarse models may produce quite different solutions, the error in the quantities of interests, $\calQ(u) - \calQ(u_0)$, where $\calQ(u)$ and $\calQ(u_0)$ are furnished by the two models, is of primary importance. 
In this regard, the following theorem, proved in \cite{oden2001goal, oden2002estimation}, relates the error $\calQ(u) - \calQ(u_0)$ in QoI to the residuals $\calR(\cdot; \cdot)$ and $\bar{\calR}(\cdot; \cdot, \cdot)$:
\begin{theorem}\label{thm:Apost}
Given any approximation $(u_0, p_0)$ of the solution $(u,p)$ of \eqref{eq:fwdAdj}, the following a-posteriori error representation holds:
\begin{equation}\label{eq:errQoITum}
\calQ(u) - \calQ(u_0) = \calR(u_0; p_0) + \frac{1}{2} \orbrack{\calR(u_0; \eps_0) + \bar{\calR}(u_0; e_0, p_0)} + r_1,
\end{equation}
where
\begin{equation*}
	\begin{split}
		r_1 &= \frac{1}{2}\int_0^1 \left\{\calQ'''(u_0 + se_0; e_0, e_0, e_0) \right. \\
		&\qquad \quad \left.- 3\calB''(u_0+se_0; e_0, e_0, \eps_0) - \calB'''(u_0+se_0; e_0, e_0, e_0, p_0 + s\eps_0)\right\} (s - 1)s \dd s .
	\end{split}
\end{equation*}
\qed
\end{theorem}

Next, Lemma 1 in \cite{oden2001goal} is called upon and slightly extended to show that higher-order approximations of $\bar{\calR}(u_0; e_0, p_0)$ in terms of $\calR(u_0; \eps_0)$ can be obtained.
\begin{lemma}\label{lem:Rbar}
Given any approximation $(u_0, p_0)$ of the solution $(u,p)$ of \eqref{eq:fwdAdj}, there holds
\begin{equation}
\bar{\calR}(u_0; e_0, p_0) = \calR(u_0; \eps_0) + r_2 \label{eq:AdjResFwdRes1}
\end{equation}
and
\begin{equation}
\bar{\calR}(u_0; e_0, p_0) = \calR(u_0; \eps_0) - \calQ''(u_0; e_0, e_0) + \calB''(u_0; e_0, e_0, p_0) + \frac{1}{2}\calB''(u_0; e_0, e_0, \eps_0) + r_3, \label{eq:AdjResFwdRes2}
\end{equation}
where
\begin{equation*}
		r_2 = \int_0^1 \left\{ \calB''(u_0 + s e_0; e_0, e_0, p_0 + s\eps_0) - \calQ''(u_0 + s e_0; e_0, e_0) \right\} \dd s
\end{equation*}
and
\begin{equation*}
		r_3 = \int_0^1 \left\{ \calB'''\left(u_0 + s e_0; e_0, e_0, e_0, p - \frac{1}{2}(1-s) \eps_0 \right) - \calQ'''(u_0 + s e_0; e_0, e_0, e_0) \right\} (1 - s) \dd s .
\end{equation*}
\end{lemma}

Equalities \eqref{eq:AdjResFwdRes1} and \eqref{eq:AdjResFwdRes2} are derived in \cref{s:proof}. Equation \eqref{eq:AdjResFwdRes1} is established in \cite{oden2001goal} and \eqref{eq:AdjResFwdRes2} is derived through straightforward algebraic manipulations described in \cref{s:proof}. 
Proof of these two equations, as shown next, leads to variational problems that can be used to compute the approximations of the errors $(e_0, \eps_0)$.

Combining \cref{thm:Apost} and \cref{lem:Rbar} to eliminate $\bar{\calR}(u_0; e_0, p_0)$, the following pair of representations of modeling error is obtained:
\begin{equation}\label{eq:errQoI1}
	\calQ(u) - \calQ(u_0) = \calR(u_0; p_0) + \calR(u_0; \eps_0) + r_1 + \frac{r_2}{2} , 
\end{equation}
and
\begin{equation}\label{eq:errQoI2}
	\calQ(u) - \calQ(u_0) = \calR(u_0; p_0) + \calR(u_0; \eps_0) - \calQ''(u_0; e_0, e_0) + \calB''(u_0; e_0, e_0, p_0  + \eps_0/2)   + r_1 + \frac{r_3}{2} .
\end{equation}
The remainder terms in the above equations, in general, depend nonlinearly on the solutions $(u_0, p_0)$ and $(u,p)$. 

While \eqref{eq:errQoI1} and \eqref{eq:errQoI2} provide exact representations of error in the QoI, we may often consider approximations of this error, as noted in \cite{oden2001goal, oden2002estimation}, that can be more easily computed. For example, $\calR(u_0; p_0)$ is readily computable whenever $u_0$ and $p_0$ are known. 
In instances in which the higher order terms, $r_1, r_2, r_3$, may be neglected when $e_0$ and $\eps_0$ are sufficiently small, the following approximate error estimators are obtained, noting that $\calR(u_0; p_0) + \calR(u_0; \eps_0) = \calR(u_0; p_0 + \eps_0) = \calR(u_0; p)$ (as $\calR(\cdot; \cdot)$ is linear in second argument and $p_0 + \eps_0 = p$):
\begin{equation}\label{eq:errQoI3}
	\calQ(u) - \calQ(u_0) \approx \calR(u_0; p) =: \Xi_1(u, u_0) 
\end{equation}
and
\begin{equation}\label{eq:errQoI4}
	\calQ(u) - \calQ(u_0) \approx \calR(u_0; p) - \calQ''(u_0; e_0, e_0) + \calB''(u_0; e_0, e_0, p_0  + \eps_0/2) =: \Xi_2(u, u_0) ,
\end{equation}
$\Xi_1(\cdot, \cdot)$ and $\Xi_2(\cdot, \cdot)$ being approximations of the right-hand side of \eqref{eq:errQoITum}. Clearly from \cref{thm:Apost} and \cref{lem:Rbar}, the remainder term in the $\Xi_1$ approximation involves $\calB''$ and $\calQ''$. On the other hand, in the $\Xi_2$ approximation, the remainder term involves $\calB''$ and $\calQ'''$.  

\subsection{Approximation of solution error}\label{ss:errApprox}
The approximations $\Xi_1$ and $\Xi_2$ of the QoI error depend on errors $(e_0, \eps_0)$ and, therefore, involve the solution $(u,p)$ of the problem \eqref{eq:fwdAdj}. 
To bypass solving the high-fidelity problem \eqref{eq:fwdAdj}, a pair of linear variational problems for approximations $(\hat{e}_0, \hat{\eps}_0)$ of $(e_0, \eps_0)$ can be obtained, following \cite{oden2001goal}. In addition, a ``second-order" variational problem for approximations $(\hat{e}_0, \hat{\eps}_0)$ of the error pair an be derived which is nonlinear (quadratic) in $\hat{e}_0$ but generally more accurate that the first-order approximation. 

Referring to the derivations given in the \cref{s:proof}, particularly \eqref{eq:RB1}, and ignoring the higher order terms, i.e. $\calB''$, there holds, for any $q\in \calV$, that $\calB'(u_0; e_0, q) \approx \calR(u_0; q)$. Similarly, ignoring the higher order terms in \eqref{eq:RbarB1}, there holds, for any $v\in \calV$,  $\calB'(u_0; v, \eps_0)  \approx \bar{\calR}(u_0; v, p_0)$. This leads to the following variational problems for $(e_0, \eps_0)$: 

Given any approximation $(u_0, p_0)$ of the solution of \eqref{eq:fwdAdj}, find $(\hat{e}_0, \hat{\eps}_0) \in \calV^2$ such that, for all $(v, q) \in \calV^2$, there holds
\begin{tbox}
	\begin{equation}\label{eq:errWeakForm1}
		\begin{split}
			\calB'(u_0; \hat{e}_0, q) &= \calR(u_0; q), \\
			\calB'(u_0; v, \hat{\eps}_0) &= \bar{\calR}(u_0; v, p_0).
		\end{split}
	\end{equation}
\end{tbox}

Similarly, referring to the derivations given in \cref{s:proof}, particularly \eqref{eq:RB2}, and ignoring the higher order term, there holds, for any $q\in \calV$, $\calB'(u_0; e_0, q) + \frac{1}{2}\calB''(u_0; e_0, e_0, q) \approx \calR(u_0; q)$. And ignoring the higher order terms in \eqref{eq:RbarB2}, for any $v\in \calV$, the following holds,
\begin{equation*}
\calB'(u_0; v, \eps_0) -\calQ''(u_0; e_0, v) + \calB''(u_0; e_0, v, p) \approx \bar{\calR}(u_0; v, p_0) 
\end{equation*}
or
\begin{equation*}
\calB'(u_0; v, \eps_0) -\calQ''(u_0; e_0, v) + \calB''(u_0; e_0, v, p_0) + \calB''(u_0; e_0, v, \eps_0) \approx \bar{\calR}(u_0; v, p_0) 
\end{equation*}
Arguing as before, the following variational problems are obtained:

Given any approximation $(u_0, p_0)$ of the solution of \eqref{eq:fwdAdj}, find $(\hat{e}_0, \hat{\eps}_0) \in \calV^2$ such that, for all $(v, q) \in \calV^2$, there holds
\begin{tbox}
	\begin{equation}\label{eq:errWeakForm2}
		\begin{split}
			\calB'(u_0; \hat{e}_0, q) + \frac{1}{2}\calB''(u_0; \hat{e}_0, \hat{e}_0, q) &= \calR(u_0; q), \\
			\calB'(u_0; v, \hat{\eps}_0) -\calQ''(u_0; \hat{e}_0, v) + \calB''(u_0; \hat{e}_0, v, \hat{\eps}_0) &= \bar{\calR}(u_0; v, p_0) - \calB''(u_0; \hat{e}_0, v, p_0).
		\end{split}
	\end{equation}
\end{tbox}

\begin{remark}
	The two versions of equations for the approximate error pair $(\hat{e}_0, \hat{\eps}_0)$, \eqref{eq:errWeakForm1} and \eqref{eq:errWeakForm2}, can be obtained directly from the forward and dual problems, \eqref{eq:fwdAdj}, by simply performing the Taylor series expansion of $\calB(u; v)$ and $\calB'(u; v, p)$ about $u_0$. For example, to obtain \eqref{eq:errWeakForm1}$_1$, subtracting $\calB(u_0; q)$ from the both sides in \eqref{eq:fwdAdj}$_1$ and proceeding as follows to get
	\begin{align*}
		& \calB(u; q) - \calB(u_0; q)  = \calF(q) - \calB(u_0; q) = \calR(u_0; q)  \\
		\Rightarrow \quad & \calB'(u_0; e_0, q) + \int_0^1 \calB''(u_0 + s e_0; e_0, e_0, q) (1-s) \dd s = \calR(u_0; q).
	\end{align*}
	In the above calculation,  \eqref{eq:errWeakForm1}$_1$ is recovered by discarding $\calB''$. The following equality is used in deriving the above relation:
	\begin{equation*}
		\calB(u; q) - \calB(u_0; q) =  \calB'(u_0; e_0, q) + \int_0^1 \calB''(u_0 + s e_0; e_0, e_0, q) (1-s) \dd s .
	\end{equation*}
	Proceeding in a similar fashion and using higher-order relations for $\calB(u; q) - \calB(u_0; q)$, \eqref{eq:errWeakForm2}$_1$ can be established. The equations for $\eps_0$ can also be established in a similar manner. 
\end{remark}

\begin{remark}
The pair of equations in \eqref{eq:errWeakForm1} define generally solvable linear variational problems for approximations of the error functions $e_0$ and $\eps_0$ that, when solved and introduced into the residuals, greatly reduce the computational cost of computing goal-oriented estimates.
\end{remark}

\begin{remark}
To compute the error in the QoIs using either the exact error representation in \eqref{eq:errQoI1} and \eqref{eq:errQoI2} or the approximate representations \eqref{eq:errQoI3} and \eqref{eq:errQoI4}, a version of the fine problem must be solved. For example, given $(u_0, p_0)$ of the coarse model, the errors $(e_0, \eps_0)$ can be computed directly by solving the fine problem \eqref{eq:fwdAdj} for $(u,p)$. 
Alternatively, \eqref{eq:errWeakForm1} can be solved to compute approximations of errors $(e_0, \eps_0)$. Out of these two choices, the latter is generally preferable as it involves decoupled linear equations for $(e_0, \eps_0)$.  
\end{remark}

\subsection{Simplified estimates using the approximation of solution error}\label{ss:errApproxSimple}
In this section, the quality of approximations $\Xi_1$ and $\Xi_2$ when the errors $(e_0, \eps_0)$ are replaced by the approximate errors $(\hat{e}_0, \hat{\eps}_0)$ is examined. By combining \cref{thm:Apost} and \cref{lem:Rbar} with the variational problem \eqref{eq:errWeakForm1}$_1$, one can show that
\begin{align}\label{eq:errQoI5}
	\calQ(u) - \calQ(u_0) &= \calR(u_0; p_0) + \frac{1}{2}R(u_0; \eps_0) + \frac{1}{2}\bar{\calR}(u_0; e_0, p_0) + r_1 \notag \\
	&= \calR(u_0; p_0) + \bar{\calR}(u_0; e_0, p_0) + r_1 - \frac{r_2}{2} \notag \\
	&= \calB'(u_0; \hat{e}_0, p_0) + \calQ'(u_0; e_0) - \calB'(u_0; e_0, p_0)  + r_1 - \frac{r_2}{2} \notag \\
	&= \calQ'(u_0; \hat{e}_0) + \left[\calQ'(u_0; e_0 - \hat{e}_0) - \calB'(u_0; e_0 - \hat{e}_0, p_0) + r_1 - \frac{r_2}{2} \right].
\end{align}
Thus, in addition to the remainder terms $r_1$ and $r_2$, additional error terms arise due to the approximation of $e_0$ by $\hat{e}_0$. Similarly, by combining \eqref{eq:errWeakForm2}$_1$ and  \cref{thm:Apost} and \cref{lem:Rbar}, the following can be shown:
\begin{align}\label{eq:errQoI6}
	\calQ(u) - \calQ(u_0) &= \calQ'(u_0; \hat{e}_0) + \frac{1}{2}\calQ''(u_0; \hat{e}_0, \hat{e}_0) \notag \\
	&\qquad + \left[\calQ'(u_0; e_0 - \hat{e}_0) + \frac{1}{2}\calQ''(u_0; e_0 + \hat{e}_0, e_0 - \hat{e}_0) - \calB'(u_0; e_0 - \hat{e}_0, p_0) \right. \notag\\ 
	&\qquad\qquad  \left.- \frac{1}{2} \calB''(u_0; e_0 + \hat{e}_0, e_0 - \hat{e}_0) + r_1 - \frac{r_3}{2} \right].
\end{align}
Effectively, in both estimates, the remainder terms (all terms inside the square brackets) depend on $\calQ'$ and $\calB'$. 

\begin{remark}
	From the above equations, it is observed that if the error approximation $\hat{e}_0$ is known and employed in \eqref{eq:errQoI5} and \eqref{eq:errQoI6}, and if the terms in the square brackets including $\calQ''(u_0; \hat{e}_0, \hat{e}_0)$ are negligible, then $\calQ'(u_0; \hat{e}_0)$ may provide a readily computable approximation of the QoI error, $\calQ(u) - \calQ(u_0)$. These additional approximations are explored in specific applications in \cref{s:examples}. 
\end{remark}

%%%%%%%%%%%%%%%%%%%%%%%%%%%%%%%%%%%%%%%%%%%%%%%%%%%%%%%%%%%%%%%%%
%%%%%%%%%%%%%%%%%%%%%%%%%%%%%%%%%%%%%%%%%%%%%%%%%%%%%%%%%%%%%%%%%
\section{Bayesian Model Calibration Using  Goal-Oriented A-Posteriori Estimates}\label{s:bayes}
The principal goals of using parameterized computational models to predict events that take place in the physical universe, as noted repeatedly, are the quantities of interests. 
In the context of the present exposition, the QoIs, $\calQ(u(\btheta))$, delivered by a ``high-fidelity model" which has parameters $\btheta \in \Theta \subset \bbR^m$, is sought. To cope with uncertainties in the observational data $y$ and the imperfection in the model itself, Baye's rule is employed and a likelihood probability density, $\pi_{like}(y | \btheta)$, describing the probability distribution of the data conditioned on the parameters $\btheta$, is sought. But the observational data to which we have access is often insufficient to reliably calibrate the high-fidelity model: it may only be used to calibrate a lower-fidelity surrogate which may deliver with acceptable accuracy a QoI, $\calQ(u_0(\btheta_0))$, $\btheta_0 \in \Theta_0 \subset \bbR^{m_0}$ being the vector of parameters in the lower-fidelity model. This low-fidelity filter of data  then provides the only apparent connection with observational data available. So, the data, $y$, is taken to be
\begin{equation}
y = \calQ(u_0(\btheta_0)) .
\end{equation}
Following standard statistical arguments \cite{oden2018adaptive, oden2017predictive, oden2017foundations}, let $g$ denote the actual physical reality of an event to be predicted by our model; i.e., the ``ground truth". Then data $y = f(g, \eps)$, $f(\cdot, \cdot)$ describing a ``noise model" and $\eps$ the experimental noise. Assuming a linear additive model, $f(g, \eps) = g + \eps$, gives
\begin{equation}\label{eq:data}
y = \calQ(u_0(\btheta_0)) = g + \eps. 
\end{equation}

The high-fidelity model predicts the truth $g$ as $\calQ(u(\btheta))$, which may differ from reality due to model inadequacy or modeling error. Assuming a linear additive models of modeling error, the following relates the model prediction to the ground truth,
\begin{equation}\label{eq:data2}
g = \calQ(u(\btheta)) + \gamma(\btheta),
\end{equation}
where $\gamma  = \gamma(\btheta)$ is the modeling error (or ``model inadequacy") which depends on the parameters $\btheta$. Combining \eqref{eq:data} and \eqref{eq:data2}, the ground truth can be eliminated to give
\begin{equation}
y - \eps = g = \calQ(u(\btheta)) + \gamma(\btheta) \qquad \Rightarrow \qquad y - \calQ(u(\btheta)) = \eps + \gamma(\btheta) =: \bar{\eps} .
\end{equation}
Noting that $y = \calQ(u_0(\btheta_0))$, the total error $\bar{\eps}$, the sum of noise and model inadequacy, is equal to the goal-oriented error. 

If $\rho (\bar{\eps})$ denotes the probability density of $\bar{\eps}$ and $\pi_{like}(y | \btheta)$ is the likelihood probability of data $y$ conditioned on given model parameters $\btheta$, then
\begin{equation*}
\pi_{like}(y | \btheta) = \rho(\bar{\eps}) = \rho(y - \calQ(u(\btheta))) = \rho(\calQ(u_0(\btheta_0)) - \calQ(u(\btheta))).
\end{equation*}
As a first approximation, it is reasonable to assume that $\rho$ is Gaussian with zero mean, $\rho \sim \calN(0, \sigma)$, so that
\begin{equation}\label{eq:likePDF}
\pi_{like}(y | \btheta) = \frac{1}{\sigma \sqrt{2\pi}} \exp\left[ -\frac{|\calQ(\btheta_0; u_0) - \calQ(\btheta; u)|^2}{2\sigma^2}\right] ,
\end{equation}
where $\sigma$ is the standard deviation. 
Thus, the likelihood depends upon the error in the QoI and is estimated using the calculations described in the previous section. 
Then the posterior probability density of the parameters $\btheta$ of the high fidelity model is given by Baye's rule,
\begin{equation}\label{eq:BayesEq}
\pi_{post}(\btheta | y) = \frac{ \pi_{like}(y | \btheta) \pi_{prior}(\btheta)}{\pi_{evid}(y)},
\end{equation}
where $\pi_{prior}(\btheta)$ is a  prior probability density of the parameters and $\pi_{evid}(y)$ is the evidence density
\begin{equation*}
\pi_{evid}(y) = \int_{\Theta} \pi_{like}(y | \btheta) \pi_{prior}(\btheta) \dd \btheta .
\end{equation*}

In computations presented in the next section, a version of MCMC methods is applied to generate samples of the posterior of the model parameters using \eqref{eq:BayesEq}.

%%%%%%%%%%%%%%%%%%%%%%%%%%%%%%%%%%%%%%%%%%%%%%%%%%%%%%%%%%%%%%%%%
%%%%%%%%%%%%%%%%%%%%%%%%%%%%%%%%%%%%%%%%%%%%%%%%%%%%%%%%%%%%%%%%%
\section{Applications}\label{s:examples}
In this section, the method of model parameter estimation described earlier is applied to two classes of problems. The first application involves a nonlinear boundary-value problem defined on a 2D domain. The model parameters of the fine (nonlinear) model are inferred using as the coarse model a simple linearized model. 
The second application is concerned with the tissue-scale tumor growth models. Specifically, a transient nonlinear partial differential equation modeling tumor growth as the fine model and a surrogate model characterized by a transient linear partial differential equation is considered. 
All computations are performed on a Macbook laptop\footnote{Further details about the architecture and software are available at this link: https://github.com/prashjha/GoalOrientedModelCalibration} (serial program execution) using the Fenics library \cite{alnaes2015fenics, logg2012automated}.

%%%%%%%%%%%%%%%%%%%%%%%%%%%%%%%%%%%%%%%%%%%%%%%%%%%%%%%%%%%%%%%%%
\subsection{Quasi-linear second order elliptic boundary-value problem}
Let $\Omega = (0,1)^2$ be an open square domain with boundary $\p \Omega$. Consider the following quasi-linear elliptic problem as the fine model: Find $u = u(\bx) \in \calV := \{v\in H^1(\Omega): u = 0 \text{ on } \p \Omega\}$ such that
\begin{equation}\label{eq:NonlinWeak}
\formBExpand{u}{v} = \int_\Omega fv \dd \bx, \qquad \forall v\in \calV,
\end{equation}
where $\kappa >0$ and $\alpha$ are parameters, $f\in \calV'$ is the source term, $\bx = (x_1, x_2) \in \Omega$, and $\dd \bx = \dd x_1 \dd x_2$. Here, it is assumed that the problem \eqref{eq:NonlinWeak} is well-posed and a unique solution exists in $\calV$. 
Associated with this problem, forms $\calB$ and $\calF$ are defined as follows
\begin{equation}
\calB(u; v) = \formBExpand{u}{v}, \qquad
\calF(v) = (f, v) = \int_\Omega fv \dd \bx .
\end{equation}
Further, the volume average of the solution $u$ is taken as the quantity of interest, i.e., 
\begin{equation}
\calQ(u) = \int_\Omega u \dd \bx .
\end{equation}
In this case, the model parameters are $\btheta = (\kappa, \alpha)$, and
\begin{equation*}
\calB'(u; v, p) = \formBDerExpand{u}{v}{p}, \qquad \calQ'(u; v) = \calQ(v).
\end{equation*}

A linear elliptic problem is taken as a coarse problem: Find $u_0 \in \calV_0 := \calV$ such that
\begin{equation}\label{eq:LinWeak}
\calB_0(u_0; v_0) = \calF(v_0), \qquad \forall v_0 \in \calV_0
\end{equation}
with
\begin{equation*}
\calB_0(u_0, v_0) = \int_\Omega \kappa_0 \nabla u_0 \cdot \nabla v_0 \dd \bx.
\end{equation*}
Here, $\kappa_0$ is the diffusion coefficient assigned a fixed value of $0.25$. The forcing function $f\in \calV'$ in both nonlinear and linear BVPs is fixed as follows
\begin{equation*}
f = f(\bx) = 10 \cos^2(4\pi x_1) \cos^2(4\pi x_2), \qquad \bx = (x_1, x_2) \in \Omega.
\end{equation*}

\subsubsection{Comparing various goal-oriented error estimates}
In this section, the various goal-oriented estimates are calculated for the example problem \eqref{eq:NonlinWeak} and compared in order to assess their accuracy. 
For the numerical approximation of PDEs in \eqref{eq:NonlinWeak} and \eqref{eq:LinWeak}, a quadrilateral finite element approximation on a mesh of $50\times 50$ elements with first-order shape functions is employed. Newton's method is used to solve the discrete nonlinear problem \eqref{eq:NonlinWeak}.  

First, the forward and dual problem \eqref{eq:fwdAdjLF} corresponding to the linear BVP \eqref{eq:LinWeak} are solved to obtain $(u_0, p_0)$. In this case, $\calQ(u_0)  = \int_\Omega u_0 \dd \bx = 0.33577$. 
Next, the QoI error estimates $\Xi_1$ and $\Xi_2$ defined in \eqref{eq:errQoI3} and \eqref{eq:errQoI4}, respectively, are compared with the ``exact" error in the QoI. To do this, the parameters in the nonlinear model \eqref{eq:NonlinWeak} are set to $\kappa = 0.25, \alpha = 10$. 
For these parameters in the nonlinear model, the compute time to solve the nonlinear problem \eqref{eq:NonlinWeak} for $u$ and the linear problem \eqref{eq:errWeakForm1} for an approximate error $\hat{e}_0$ are found to be $1.24$ and $0.557$ seconds, respectively.

\begin{itemize}
\item[\textbf{(i)}] \textbf{Calculating estimates using the solution of the fine model.}  In this case, $(u,p)$ is obtained by solving \eqref{eq:fwdAdj} for the example corresponding to \eqref{eq:NonlinWeak}. Using $u$, $\calQ(u)$ and the ``exact" QoI error are found to be
\begin{equation*}
\calQ(u) = 0.1163 \quad \Rightarrow \quad \calQ(u) - \calQ(u_0) = \calQ(e_0) = 0.1163  - 0.33577  = -0.21947,
\end{equation*}
where ``exact" refers to the fact that $\calQ(u)$ is evaluated using the finite element solution $u$ of the high-fidelity model; ``exact" is inside the double quote to note that there is still a numerical discretization error in $\calQ(u)$ and in $\calQ(u) - \calQ(u_0)$. 

With $(u,p)$ in hand and $(u_0, p_0)$ known, the errors are readily available and using these the two estimates of the QoI error are found to be
\begin{equation*}
	\begin{split}
		\calQ(u) - \calQ(u_0) \approx \Xi_1(u, u_0) &= -0.2069, \\
		\calQ(u) - \calQ(u_0) \approx \Xi_2(u, u_0) &= -0.22468 .
	\end{split}
\end{equation*}

\item[\textbf{(ii)}] \textbf{Calculating estimates using the approximate errors.} In this case, the approximate errors $(\hat{e}_0, \hat{\eps}_0)$ are computed by solving \eqref{eq:errWeakForm1}. Using $(\hat{e}_0, \hat{\eps}_0)$, $u$ and $p$ can be approximated as $u \approx u_0 + \hat{e}_0$ and $p \approx p_0 + \hat{\eps}_0$.  The QoI and the error in the QoI are
\begin{equation*}
\calQ(u) = 0.12306 \quad \Rightarrow \quad \calQ(u) - \calQ(u_0) = \calQ(\hat{e}_0) = 0.12306  - 0.33577  = -0.21271 .
\end{equation*}

On the other hand, the estimates of the QoI errors are
\begin{equation*}
	\begin{split}
		\calQ(u) - \calQ(u_0) \approx \Xi_1(u, u_0) &= -0.21272, \\
		\calQ(u) - \calQ(u_0) \approx \Xi_2(u, u_0) &= -0.21272 .
	\end{split}
\end{equation*}
Clearly, the two estimates are very close to the ``exact" error, $\calQ(u) - \calQ(u_0) = -0.21947$. 
\end{itemize}

In the above results, approach \textbf{(ii)}, in which $(\hat{e}_0, \hat{\eps}_0)$ are utilized, is clearly more efficient than \textbf{(i)} as in this case only the linear problem \eqref{eq:errWeakForm1} is solved for errors $(\hat{e}_0,  \hat{\eps}_0)$. For the case in which errors $(\hat{e}_0, \hat{\eps}_0)$ are used to estimate the QoI error, the three estimates, $\calQ(\hat{e}_0)$, $\Xi_1$, and $\Xi_2$, show good agreement with the exact QoI error of $-0.21947$; however, computing $\calQ(\hat{e}_0)$ is simpler and more straightforward and requires solving for only the forward error, $\hat{e}_0$. Therefore, in the calibration calculations below, $\calQ(\hat{e}_0)$, is used as the estimate of the QoI error and only variational problem \eqref{eq:errWeakForm1}$_1$ is solved, yielding $\hat{e}_0$. 

\subsubsection{Calibration of the fine model} \label{sss:app1Calib}
Let $\btheta_0 = (\kappa_0)$ and $\btheta = (\kappa, \alpha)$. A log-normal prior $\pi_{prior}(\btheta)$ for $\btheta$ with ln-mean and ln-std of $(-0.6535, 2.5475)$ and $(0.1997, 0.5003)$, respectively, is assumed and the standard deviation of noise, $\sigma$, is set to $0.01$. 
With the likelihood defined in \eqref{eq:likePDF}, the Hippylib \cite{villa2021hippylib, villa2018hippylib} library is used to compute the approximate posterior probability density function $\pi_{post}(\btheta | y)$ using the Bayes rule. To generate samples, four MCMC chains are used with the maximum number of samples drawn for each chain set to $5000$. The total number of posterior samples from all chains after discarding $50$ percent of initial accepted samples (burn-in) is $2766$. 

\begin{figure}[ht]
     \centering
    \begin{subfigure}[t]{0.6\textwidth}
        \raisebox{-\height}{\includegraphics[width=\textwidth]{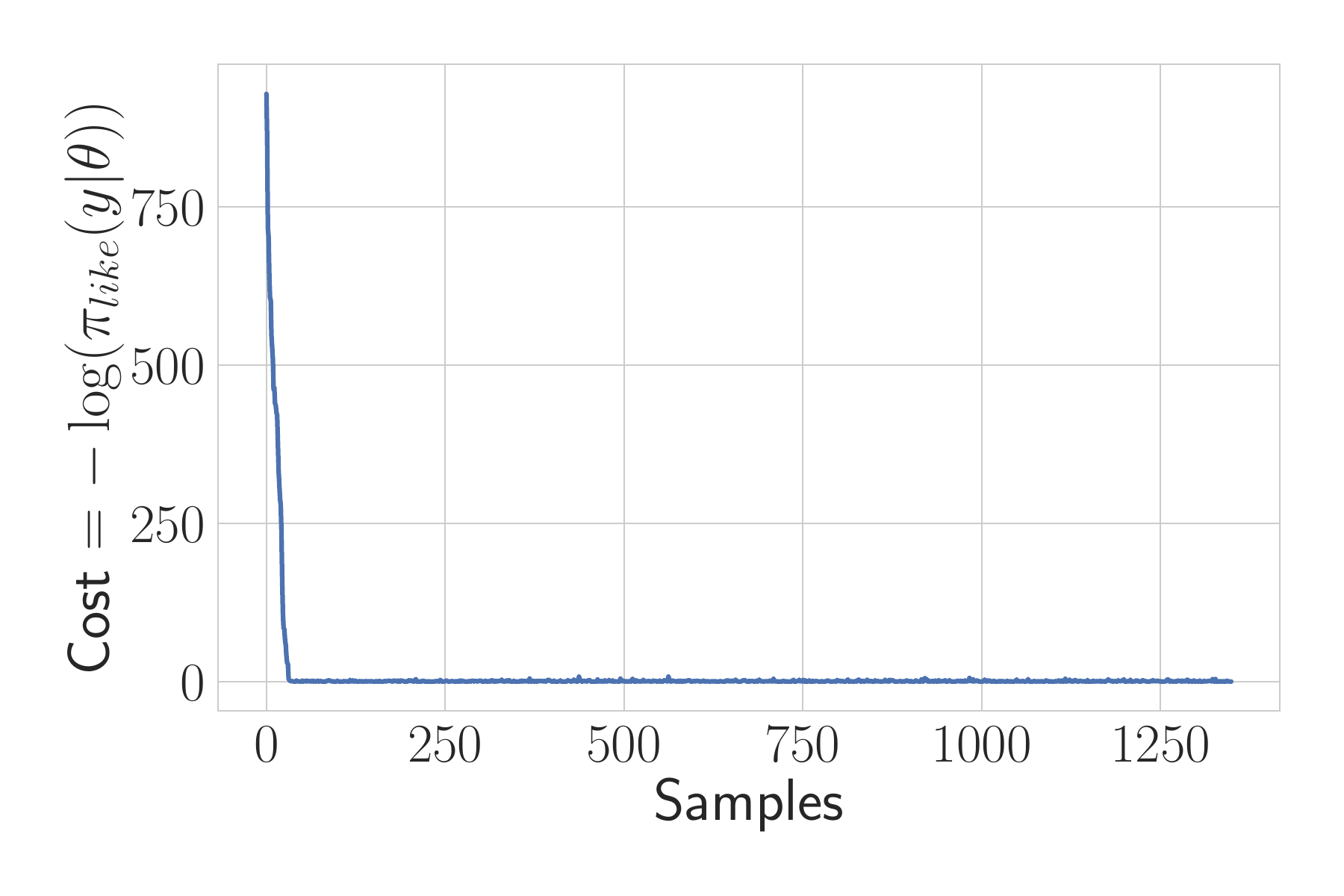}}
		\vspace{-10pt}
        \caption{Cost associated to the accepted samples.}
    \end{subfigure}
    \\
    \begin{subfigure}[t]{0.45\textwidth}
        \raisebox{-\height}{\includegraphics[width=\textwidth]{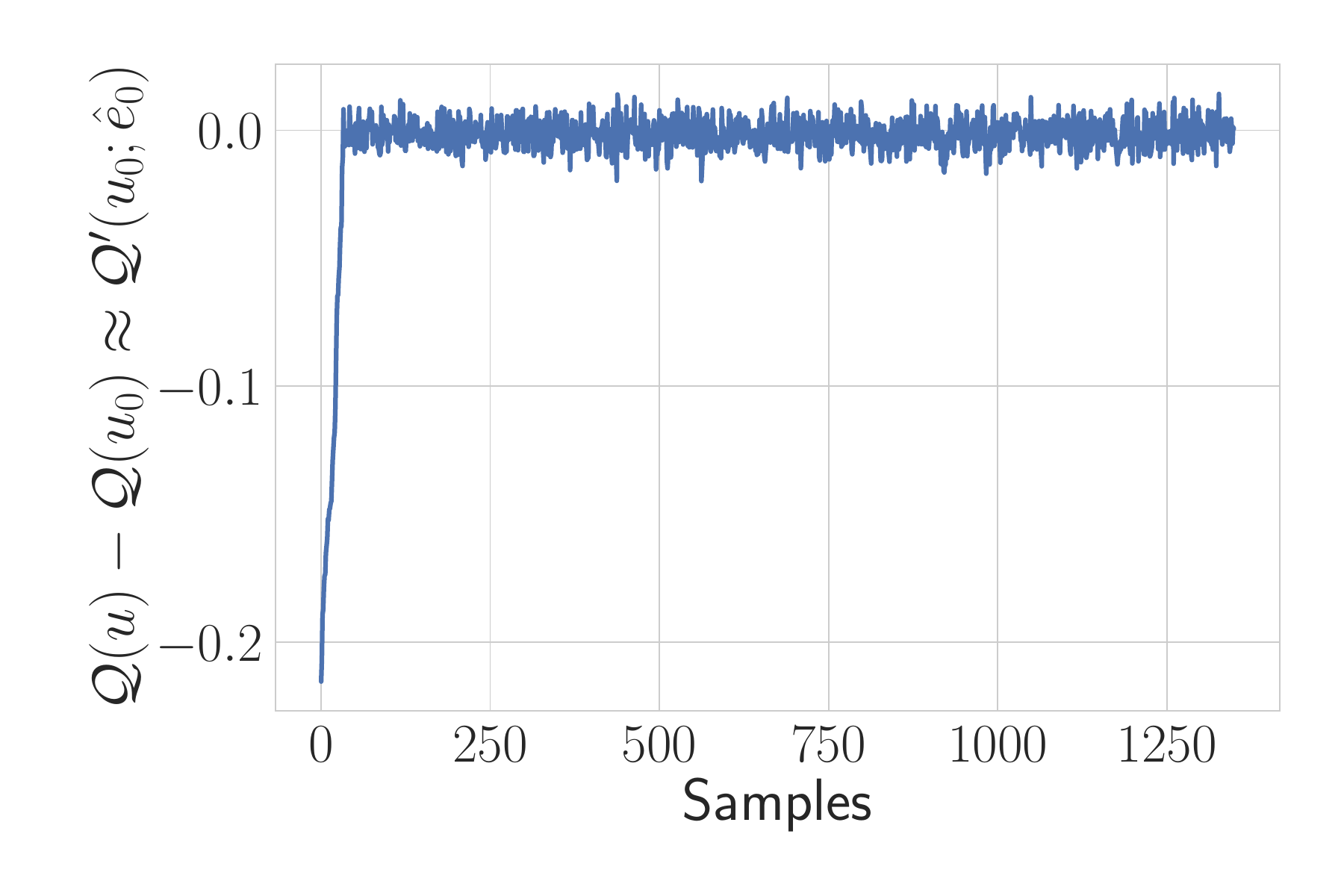}}
		\vspace{-5pt}
        \caption{QoI error $\calQ(u) - \calQ(u_0)$ using the estimate $\calQ'(u_0; \hat{e}_0)$.}
    \end{subfigure}
	\hfill
    \begin{subfigure}[t]{0.45\textwidth}
        \raisebox{-\height}{\includegraphics[width=\textwidth]{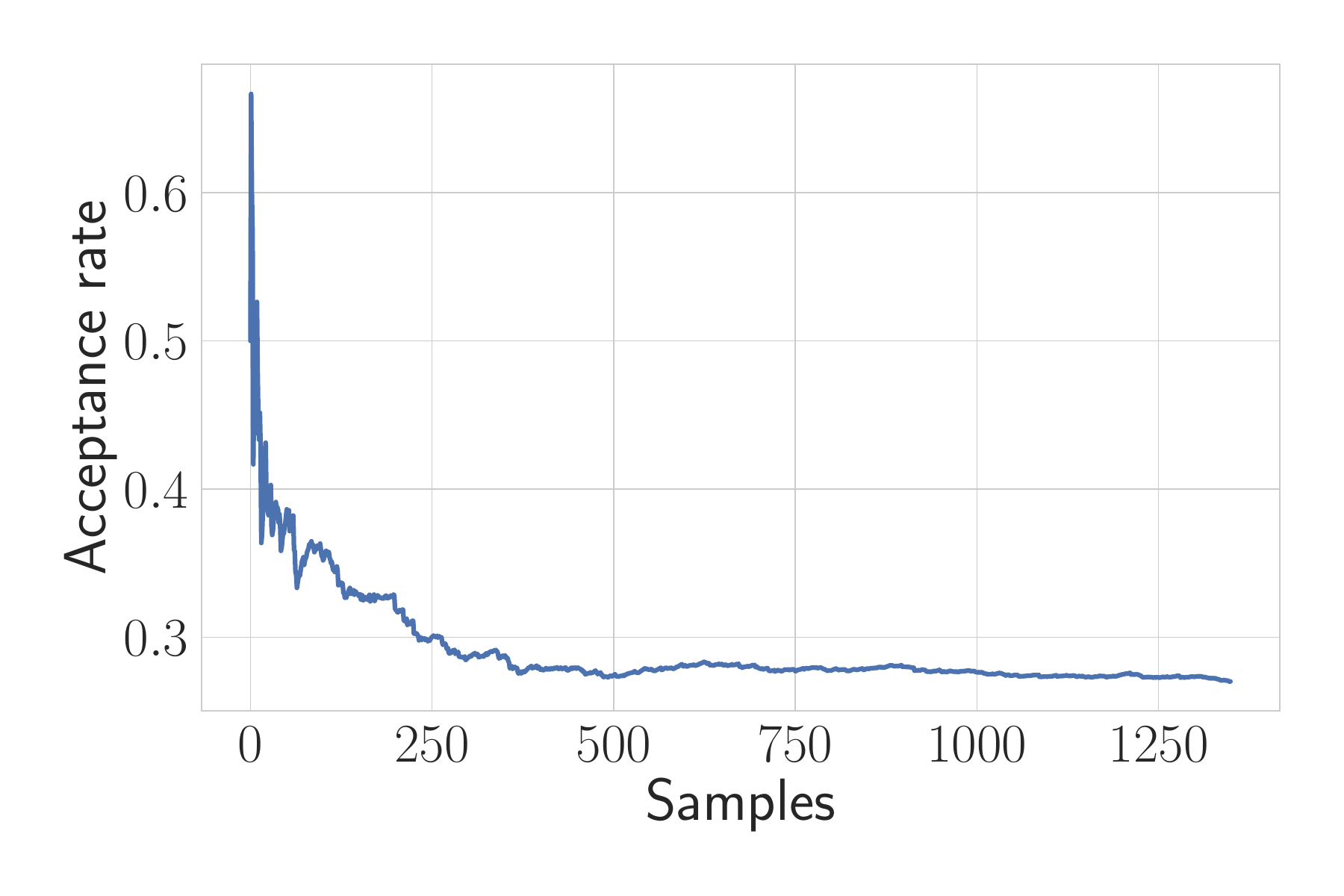}}
		\vspace{-5pt}
        \caption{Sample acceptance rate.}
    \end{subfigure}
    \caption{Results from one MCMC chain. In (A), the variation in cost associated with the accepted samples is shown. In (B), the QoI error associated with the accepted samples is shown. The rate of sample acceptance is shown in (C). Initially, the cost and the error in QoI is very high, but the rate quickly stabilizes so that the mean approaches a constant value.}\label{fig:mcmcStats}
\end{figure}

\begin{figure}[ht]
     \centering
    \includegraphics[width=\textwidth]{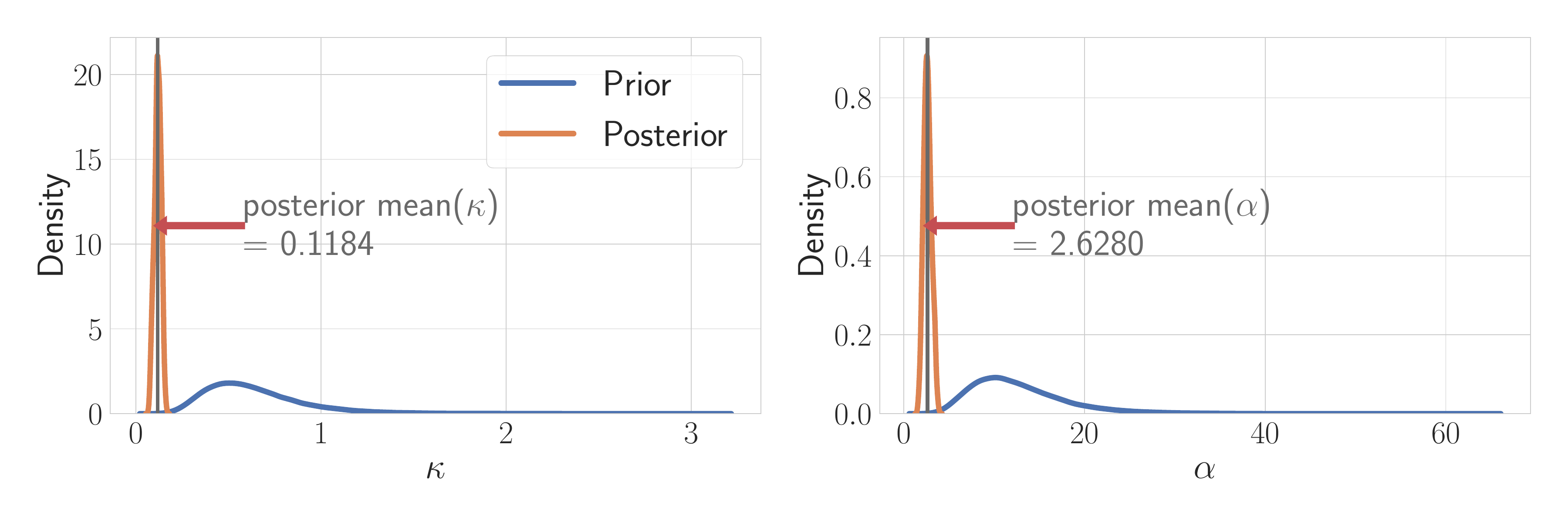}
	\vspace{-20pt}
    \caption{The prior and posterior probability density functions and mean values of posterior samples for parameters in the high-fidelity model. }\label{fig:comparePriorPost}
\end{figure}

In \cref{fig:mcmcStats}, the cost, error in QoI ($\calQ(\btheta, u) - \calQ(\btheta_0, u_0) \approx \calQ'(u_0; \hat{e}_0) = \calQ(\hat{e}_0)$), and the sample acceptance rate during the MCMC computation for one of the chains are shown. 
In \cref{fig:comparePriorPost}, the prior and posterior densities are compared. The mean and the standard deviation of the posterior samples are $\mu_{\btheta} = (0.118, 2.628)$ and $\sigma_{\btheta} = (0.018, 0.433)$, respectively. Particularly, the QoI with parameters $\btheta = \mu_{\btheta}$ is $\calQ(u) = 0.334$ and the error in QoI, $\calQ(u) - \calQ(u_0)$, is $-0.0013$ ($-0.4$ percent of $\calQ(u_0)$). At the mean parameter $\mu_{\btheta}$, the approximate QoI error is $\calQ(\hat{e}_0) = -0.0016$, i.e., $-0.46$ percent of $\calQ(u_0)$. 

\subsubsection{Reliability of the calibration under the QoI error approximation}
In the calibration of the model in the previous section, the estimate $\calQ(\hat{e}_0)$ of $\calQ(u) - \calQ(u_0)$ is used and, therefore, solving the fine problem was avoided. 
Obviously, the use of such an approximation in the calibration steps as opposed to the ``exact" in the QoI could affect the accuracy of the posterior samples. This leads to the question: How do the two posterior samples, one in which the approximate QoI error is employed and another in which the QoI error is computed exactly (up to discretization error), differ?

\begin{figure}[ht]
	\centering
	\includegraphics[width=\textwidth]{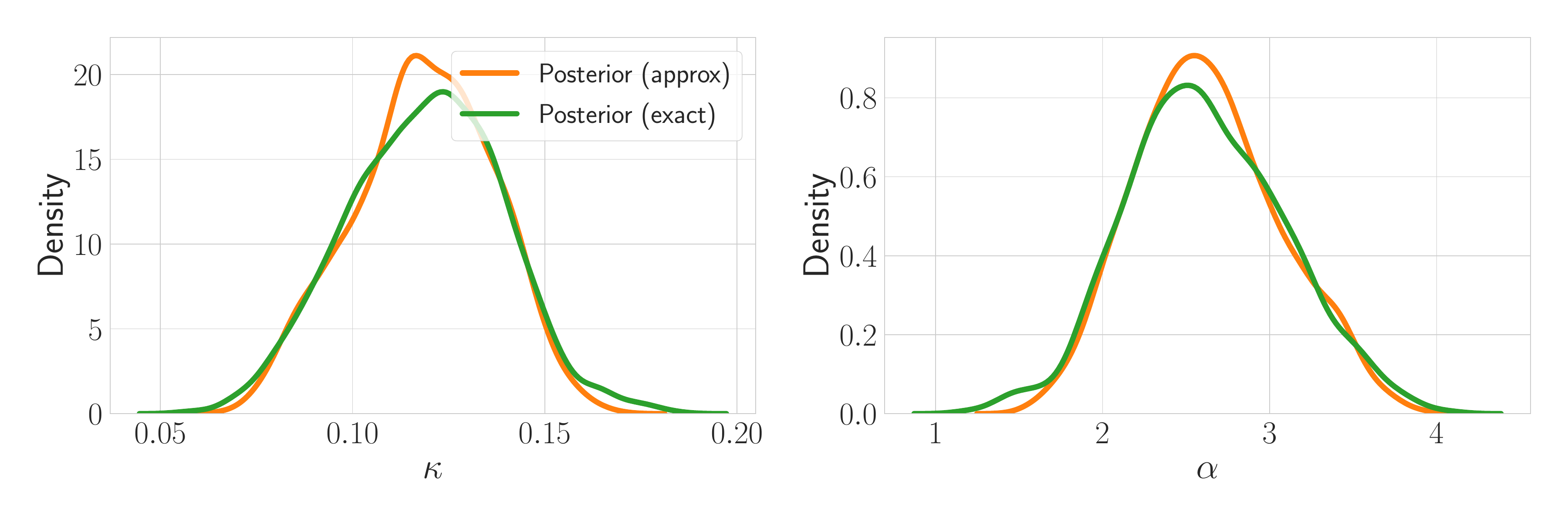}
	\vspace{-20pt}
	\caption{The posterior probability density functions for the case in which the QoI error is approximated by $\calQ(\hat{e}_0)$ and in which the QoI error is computed exactly by solving the fine problem.}\label{fig:comparePost}
\end{figure}

To answer this, the MCMC sampling described in \cref{sss:app1Calib} is repeated but now the QoI difference $\calQ(u) - \calQ(u_0)$ is computed exactly by solving the fine problem for $u$. The posterior densities for the two cases are shown in \cref{fig:comparePost}. The mean and standard deviation of posterior samples in two cases are 
\begin{align*}
	&\textbf{(approx)} \qquad &  & \mu_{\btheta} = (0.118, 2.628), \qquad \sigma_{\btheta} = (0.018, 0.433), \notag \\
	&\textbf{(``exact")} \qquad &  & \mu_{\btheta} = (0.119, 2.616), \qquad \sigma_{\btheta} = (0.02, 0.475).
\end{align*}
In this example, the estimate $\calQ(\hat{e}_0)$ of the QoI error is seen to produce a reasonably accurate representation of the posteriors for the parameters $\btheta$. 
The compute time to finish the first chain (5000 samples) in MCMC simulations in the two cases of ``exact" and approximate involving QoI error estimate are $17,760$ and $9,300$ seconds, respectively. Thus, the case in which the QoI error is approximated using goal-oriented estimates is about $1.9$ times faster than that in which the ``exact" solution is employed.

%%%%%%%%%%%%%%%%%%%%%%%%%%%%%%%%%%%%%%%%%%%%%%%%%%%%%%%%%%%%%%%%%
\subsection{Tumor growth model}
An application of these modeling error estimation methods of particular interest in the present study is that of modeling tumor growth in living tissue, an application in which data for calibration of model parameters is notoriously difficult to obtain. 
Towards this goal, a high-fidelity model based on a nonlinear PDE, the so-called Allen-Cahn phase field equation, introduced originally in the context of modeling phase changes in binary alloys \cite{allen1979microscopic}, is first presented. Phase-field models including Allen-Cahn and Cahn-Hilliard equations have been used extensively to model the tissue-scale tumor growth, see \cite{lorenzo2016tissue, alsayed2022optimal, fritz2021analysis, hawkins2012numerical, oden2010general, lima2016selection, lima2017selection, fritz2021modeling}. 
The model in this application describes the evolution of tumor volume with volume fraction, $u = u(t, \bx)$, where $t\in [0, t_F]$ is the time interval of interest and $\bx \in \Omega \subset \bbR^2$ is the point on the tissue domain $\Omega$. 
Next, a low-fidelity model based on a linear PDE describing the evolution of tumor volume fraction $u_0 = u_0(t, \bx)$ is considered. This low-fidelity model is based on the assumption that tumor growth is governed by diffusion and with proliferation defined as a linear function of the tumor volume fraction. 

%-------------------------------------------------------------------------%
\subsubsection{Models of tumor growth}
Consider a two-dimensional tissue-scale model of tumor growth over a time interval $(0, t_F)$ in which a colony of tumor cells occupies a domain $\Omega = (0,1)^2$. 
The tumor volume fraction $u = u(t, \bx)$ is assumed to be governed by the Allen-Cahn equation, for all $t \in (0, t_F), \bx \in \Omega$, 
\begin{equation}\label{eq:HFMassBal}
    \frac{\partial u}{\partial t} = \epsilon \nabla \cdot (\nabla u) - \Psi'(u) + \lambda^p u (1 - u) f - \lambda^d u,
\end{equation}
where $\epsilon$ is the interfacial width between the tumor and the external cellular tissue, $\Psi = \Psi(u) = C u^2 (1-u)^2$ is a double-well energy, $\lambda^p$ the proliferation rate of tumor cells, $f = f(t, \bx)$ a nutrient source, and $\lambda^d$ the death rate of tumor cells. In this application, a fixed form of $f$ is prescribed as
\begin{equation*}
	f(t, \bx) = \exp[-1.5 x_1], \qquad \bx = (x_1, x_2) .
\end{equation*}
In more general cases, $f$ could represent the solution of a reaction-diffusion equation governing the evolution of nutrients in the tissue domain \cite{fritz2021analysis, hawkins2012numerical, oden2010general, lima2016selection, lima2017selection, fritz2021modeling}. A homogeneous Neumann boundary condition is assumed for $u$, i.e., 
\begin{equation*}
    \nabla u \cdot \bn = 0, \qquad \forall \bx \in \partial \Omega, \forall t,
\end{equation*}
where $\bn$ is the unit outward normal to the boundary $\partial \Omega$. At $t=0$, the tissue is assumed to carry a spherical tumor of radius $r_c$; the initial condition for $u$ is taken as
\begin{equation*}
    u(0, \bx) = \bar{u}(\bx) := \begin{cases}
    	1, \qquad \text{ if } |\bx - \bx_c| < r_c, \\
    	0, \qquad \text{ otherwise},
    \end{cases}
\end{equation*}
for $\bx \in \Omega$ with $\bx_c = (0.5, 0.5)$, $r_c = 0.2821$ (so that $\int_\Omega \bar{u}(\bx) \dd\bx = 0.25$). 

\paragraph{Low-fidelity model}
Suppose $u_0 = u_0(t, \bx)$ is the tumor volume fraction obtained by solving the linear reaction-diffusion equation given by, 
\begin{equation}\label{eq:LFMassBal}
	\frac{\partial u_0}{\partial t} = D \nabla \cdot (\nabla u_0) + \lambda^p_0 u_0 f - \lambda^d_0 u_0, \qquad \forall t \in (0, t_F], \, \bx \in \Omega, 
\end{equation}
where $D$ is the diffusivity of tumor cells, $\lambda^p_0$ the proliferation rate, and $\lambda^d_0$ the death rate. It is assumed that $u_0$ satisfies the same boundary and initial condition as $u$. 

%-------------------------------------------------------------------------%
\subsubsection{Weak formulation}
To cast the PDE-based model into an appropriate functional setting, the function spaces $\calU$ and $\calV$ are introduced as follows
\begin{equation}
\calU = L^2(0, t_F; H^1(\Omega)), \qquad \calV = \{v\in \calU: \p_t v \in \calU'\},
\end{equation}
where $\calU' =  L^2(0,t_F; H^1(\Omega)')$ is the dual of $\calU$ (similarly, $H^1(\Omega)'$ is the dual of $H^1(\Omega)$). The norm of $v\in \calU$ is given by
\begin{equation*}
    \normX{v}{\calU}^2 = \int_0^{t_F} \normX{v(t)}{H^1(\Omega)}^2 \dd t, \qquad \normX{u}{H^1(\Omega)}^2 = \normX{u}{L^2(\Omega)}^2 + \normX{\nabla u}{L^2(\Omega)}^2
\end{equation*}
and 
\begin{equation*}
    \normX{w}{\calV}^2 = \normX{w}{\calU}^2 + \normX{\p_t w}{\calU'}^2,
\end{equation*}
where the norm of $v\in \calU'$ is
\begin{equation*}
    \normX{v}{\calU'}^2 = \int_0^{t_F} \normX{v(t)}{H^1(\Omega)'} \dd t, \qquad \normX{v}{H^1(\Omega)'} = \sup_{w\in H^1(\Omega)} \frac{<v, w>}{\normX{w}{H^1(\Omega)}}.
\end{equation*}
Here $<v, w>$ denotes duality pairing on $H^1(\Omega)' \times H^1(\Omega)$. 

Let $\calV_0 = \calV$, the solution space for the low-fidelity model, and let us assume that the weak solutions $u$ and $u_0$ of \eqref{eq:HFMassBal} and \eqref{eq:LFMassBal}, respectively, are in $\calV$ and $\calV_0$. 
The semilinear form $\calB : \calV \times \calV \to \bbR$ associated with \eqref{eq:HFMassBal} is given by
\begin{equation}\label{eq:formB}
    \calB(u; v) = \weakDot{u(0)}{v(0)} + \int_0^{t_F} \dualDot{\p_t u}{v}\dd t + \calA(u, v) + \calN(u; v),
\end{equation}
where $\calA: \calV \times \calV \to \bbR$ is the bilinear form and $\calN: \calV \times \calV \to \bbR$ is the semilinear form defined as
\begin{equation*}
    \calA(u, v) = \int_0^{t_F} \left\{ \weakDot{\epsilon \nabla u}{\nabla v} + \weakDot{\lambda^d u}{v}  \right\} \dd t
\end{equation*}
and
\begin{equation*}
    \calN(u;v) = \int_0^{t_F} \leftcr  - \weakDot{\lambda^p u (1 - u) f}{v} + \weakDot{\Psi'(u)}{v} \rightcr \dd t .
\end{equation*}
The linear form $\calF: \calV \to \bbR$ is taken to be
\begin{equation}\label{eq:formF}
    \calF(v) = \weakDot{\bar{u}}{v(0)} .
\end{equation}

With the above notations, the weak form of \eqref{eq:HFMassBal} reduces to the general form given by \eqref{eq:fwdAdj}$_1$. The parameters in the high-fidelity model are $\btheta = (\lambda^p, \lambda^d, \epsilon, C)$.
A bilinear form $\calB_0 : \calV \times \calV \to \bbR$ associated with the low-fidelity problem \eqref{eq:LFMassBal} is defined as, for $u_0, v_0 \in \calV_0$,
\begin{equation*}
\calB_0(u_0, v_0) =  \weakDot{u_0(0)}{v_0(0)} + \int_0^{t_F} \left\{ \dualDot{\p_t u_0}{v_0} + \weakDot{D \nabla u_0}{\nabla v_0} - \weakDot{\lambda^p_0 u_0 f}{v_0} + \weakDot{\lambda^d_0 u_0}{v_0}  \right\} \dd t.
\end{equation*}
With this $\calB_0$ and $\calF$, the weak form of \eqref{eq:LFMassBal} is given by \eqref{eq:fwdAdjLF}$_1$. The parameters in the linear model are $\btheta_0 = (\lambda^p_0, \lambda^d_0, D)$. 

For the adjoint formulation, and to compute errors using \eqref{eq:errWeakForm1} or \eqref{eq:errWeakForm2}, $\calB'(u;v, p)$ and $\calB''(u: q, v, p)$ are needed. It can be shown that, for $u, v, p \in \calV$,
\begin{equation}\label{eq:derFormB}
	\begin{split}
		    \calB'(u; v, p) &= \lim_{\eta \to 0} \eta^{-1} \left[\calB(u+\eta v; p) - \calB(u; p) \right] \\
		    & = \weakDot{p(t_F)}{v(t_F)} - \int_0^{t_F} \dualDot{\p_t p}{v} \dd t + \calA(v, p) + \calN'(u; v, p),
	\end{split}
\end{equation}
where
\begin{equation*}
   \calN'(u; v, p) = \int_0^{t_F} \left\{ \weakDot{\Psi''(u)v}{p} 
			- \weakDot{\lambda^p (1 - 2 u) v f}{p} \right\} \dd t.
\end{equation*}
Further, for any $u, v, p, q\in \calV$, it can be shown that
\begin{equation*}
	\begin{split}
		\calB''(u; q, v, p) &= \lim_{\eta \to 0} \eta^{-1} \left[\calB'(u+\eta q; v, p) - \calB(u; v, p) \right] = \calN''(u; q, v, p) \\
		&=	 \int_0^{t_F} \left\{ \weakDot{\Psi'''(u)vq}{p} 
		+ \weakDot{2 \lambda^p vqf}{p} \right\} \dd t.
	\end{split}
\end{equation*}

To establish \eqref{eq:derFormB}, the following identity, given in \cite{van2011goal}, is essential: for any function $u\in \calV$, it is true that $u \in C([0,t_F]; L^2(\Omega))$ and, therefore, the following integration by parts formula holds, for any $u, v \in \calV$,
\begin{equation*}
    \int_0^{t_F} \dualDot{\p_t u}{v} \dd t = \weakDot{u(t_F)}{v(t_F)} - \weakDot{u(0)}{v(0)} - \int_0^{t_F} \dualDot{\p_t v}{u} \dd t .
\end{equation*}

%-------------------------------------------------------------------------%
\subsubsection{Quantity of interests and the adjoint formulation}\label{sss:TumQoI}
Consider a general linear quantity of interest functional $\calQ: \calV \to \bbR$ of the form
\begin{align}\label{eq:formQ}
    \calQ(u) &= \weakDot{\bar{q}}{u(t_F)} + \int_0^{t_F} \dualDot{\tilde{q}}{u} \dd t,
\end{align}
where $\bar{q} \in H^1(\Omega)$ and $\tilde{q} \in \calU'$ are given a-priori. In the above equation, $\weakDot{\cdot}{\cdot}$ denotes the $L^2(\Omega)$-inner product and $\dualDot{\cdot}{\cdot}$ duality-pairing on $\calU'\times \calU$. 

In the definition of $\calQ$, $\bar{q}$ and $\tilde{q}$ are given fixed functions that characterize the tumor volume average at the final time and the temporal average of tumor volume averages at selected times as the QoI, 
\begin{equation}\label{eq:formQTum}
\bar{q}(\bx) = \frac{1}{|\Omega|} \qquad \text{ and } \qquad \tilde{q}(t, \bx) = \frac{1}{|\Omega|} \sum_{i=1}^{N_{o}} \frac{1}{\Delta \tau_{o}} \chi_{[\tau_{o, i}, \tau_{o, i} + \Delta \tau_{o}]}(t) ,
\end{equation}
where $|\Omega| = meas(\Omega)$, $N_o$ the number of observation time points, $\Delta \tau_o$ the temporal width to compute time average, $\chi_A = \chi_A(t)$ is the indicator function of set $A$ such that $\chi_A(t) = 1$ if $t\in A$ and $0$ otherwise, and $0\leq \tau_{o,i} \leq t_F - \Delta \tau_o$, $i=1, 2, ..., N_o$, observation times. With above definitions, for $u \in \calV$, the term
\begin{equation*}
\int_0^{t_F} \dualDot{\tilde{q}}{u} \dd t = \sum_{i=1}^{N_o} \frac{1}{\Delta \tau_o} \int_{\tau_{o, i}}^{\tau_{o,i} + \Delta \tau_{o}} \left[\frac{1}{|\Omega|} \int_\Omega u(t, \bx) \dd \bx \right] \dd t 
\end{equation*}
is the sum of the temporal average of volume average of $u$ at observation times $\tau_{o,i}$.

%-------------------------------------------------------------------------%
\paragraph{The adjoint problem} Following the standard procedure described in \cref{s:goal}, the adjoint problem reads:

Given the solution $u \in \calV$ the of forward problem \eqref{eq:fwdAdj}$_1$ associated with the high-fidelity tumor model \eqref{eq:HFMassBal}, find $p \in \calV$ such that for all $v\in \calV$
\begin{equation*}
\calB'(u; v, p) = \calQ'(u; v) = \calQ(v)
\end{equation*}
or, in expanded form,
\begin{align}
& \weakDot{p(t_F)}{v(t_F)} -\int_0^{t_F} \dualDot{\p_t p}{v} \dd t + \calA(v,p) + \calN'(u; v, p) = \weakDot{\bar{q}}{v(t_F)} + \int_0^{t_F} \dualDot{\tilde{q}}{v} \dd t .
\end{align}

This completes the derivation of a weak formulation of the high-fidelity tumor model.

%-------------------------------------------------------------------------%
\subsection{Verification of the estimates and calibration of the tumor growth model}
The parameters in the low-fidelity tumor model are assumed to be known with reasonable accuracy. To demonstrate the proposed model calibration approach, a non-dimensional setting with the tissue domain $\Omega = (0,1)^2$ and time $t$ in the interval $(0, 1)$ is considered. For the spatial discretization, a uniform mesh of quadrilateral elements of size $h = 1/50$ is considered and continuous first order finite element approximations on this mesh are employed. To discretize the problem in time, a first-order semi-implicit time marching scheme with the time step of $\Delta t = 0.005$ is employed for the nonlinear PDE whereas for the linear problem to compute $\hat{e}_0$ a first-order implicit time marching scheme with the same time step of $\Delta t = 0.005$ is used. For the nonlinear transient PDE, to solve the nonlinear problem at every time step, Picard iteration is employed. Further details on the numerical discretization are given in \cref{s:numdisc}. The known values of parameters in the low-fidelity model and parameters associated with the high-fidelity model are listed in \cref{tab:LFparams}.
\begin{table}[h!]
  \begin{center}
  \caption{Model parameters for the low and high fidelity models. For the high-fidelity model parameters, $\btheta = (\lambda^p, \lambda^d, \epsilon, C)$, a log-normal prior, i.e., $\ln(\btheta) \sim \calN(\mu, \Sigma^2)$, is assumed, where $\mu = \ln(\bar{\btheta}) +  (0.16, 0.16, 0.16, 0.16)$, $\bar{\btheta} = (0.5, 0.1, 0.01, 1)$, and $\Sigma^2 = \mathrm{diag}(0.16, 0.16, 0.16, 0.16)$.}
  \label{tab:LFparams}
    \begin{tabular}{|c|c|c|} 
      \hline
      \textbf{Parameter} & \textbf{Value} & \textbf{Description} \\
      \hline
      \hline
      $\lambda^p_0$ & 0.2 & Growth rate (LF) \\
      $\lambda^d_0$ & 0.1 & Death rate (LF) \\
      $D$ & 0.05 & Diffusivity (LF)  \\
      $\lambda^p$ & $\sim \mathrm{Lognormal}$ & Growth rate (HF) \\ 
      $\lambda^d$ & $\sim \mathrm{Lognormal}$ & Death rate (HF) \\ 
      $\epsilon$ & $\sim \mathrm{Lognormal}$ & Interfacial energy constant (HF) \\
      $C$ & $\sim \mathrm{Lognormal}$ & Double-well energy constant (HF) \\
      $N_{o}$ & 4 & Number of QoI observation time points \\
	  $\{\tau_{o, i}\}_{i=1}^{N_o}$ & $\{0.2, 0.4, 0.6, 0.8\}$ & Time points in QoI \\
	  $\Delta \tau_o$ & $0.05$ & Time interval in QoI \\
      \hline
    \end{tabular}
  \end{center}
\end{table}

\subsubsection{Goal-oriented error estimates} As in the first application, the accuracy of goal-oriented estimates, specifically $\calQ'(u_0; \hat{e}_0)$,  is first verified. Towards this end, the parameters in the high-fidelity model are assigned specific values; let $\btheta = \btheta_{test} = (0.5, 0.1, 0.01, 1)$. The low-fidelity parameters $\btheta_0 = (\lambda^p_0, \lambda^d_0, D)$ and remaining parameters are fixed and given in \cref{tab:LFparams}. At the outset, the low-fidelity problem is solved for $(u_0, p_0)$ and the QoI from the low-fidelity model is computed as $\calQ(u_0) = 1.143$.

\paragraph{Exact QoI error (up to discretization error)} The high-fidelity problem with the parameters $\btheta_{test}$ is solved and the QoI is found to be $\calQ(u) = 1.059$. And the error in the QoI is $\calQ(u) - \calQ(u_0) = -0.084$. Next, the approximate error, $\hat{e}_0$, is computed by solving \eqref{eq:errWeakForm1}$_1$. With $\hat{e}_0$ in hand, the QoI error is estimated to be
\begin{equation*}
	\calQ(u) - \calQ(u_0) \approx \calQ'(u_0; \hat{e}_0) = \calQ(\hat{e}_0) = -0.097,
\end{equation*}
and, since $u \approx u_0 + \hat{e}_0$, $\calQ(u) \approx 1.046$. 

Thus, use of the approximate error $\hat{e}_0$ produces a QoI agreeing with that of the HF model up to two digits of accuracy (difference between the approximate $\calQ(u)$ and the exact, ignoring discretization error, $\calQ(u)$ is within $1.3$ percent of the ``exact" $\calQ(u)$). Further, the approximate value of $-0.097$ for $\calQ(u) - \calQ(u_0)$ differs by only $1.2$ percent of the ``exact" $\calQ(u)$. 
These results encourage us to consider the computationally cheaper problem of solving for approximate error $\hat{e}_0$ and using it to approximate $\calQ(u) - \calQ(u_0)$. This route is followed in the model calibration results presented in the next section. 
It is noted that the compute time to solve the nonlinear problem for $u$ and computing the QoI functional $\calQ$ was found to be $403.67$ seconds while the compute time to solve the linear problem for an approximate error $\hat{e}_0$ and the QoI error estimate $\calQ(\hat{e}_0)$ was found to be $96.57$ seconds. Therefore, the inference using the QoI estimate is expected to speed up the MCMC sampling by a factor of about $4$. 

%-------------------------------------------------------------------------%
\subsubsection{Calibration of the high-fidelity tumor model}
Let $\btheta_0$ be fixed according to \cref{tab:LFparams}. A log-normal prior for the high-fidelity model parameters $\btheta = (\lambda^p, \lambda^d, \epsilon, C)$, i.e., $\ln(\btheta) \sim \calN(\mu, \Sigma^2)$, is considered, where $\mu = \ln(\bar{\btheta}) +  (0.16, 0.16, 0.16, 0.16)$, $\bar{\btheta} = (0.5, 0.1, 0.01, 1)$, and $\Sigma^2 = \mathrm{diag}(0.16, 0.16, 0.16, 0.16)$. For the standard deviation of noise, let $\sigma_{noise} = 0.01$. 
As in the case of the first application, the HippyLib \cite{villa2021hippylib, villa2018hippylib} library is used for the MCMC simulation. 
Four MCMC chains with the maximum number of the samples drawn in each chain fixed to $5000$ are considered. 
After discarding the initial $50$ percent of the samples (burn-in) in each of the chain, a total of $1272$ posterior samples were obtained from the MCMC simulations.

\begin{figure}[ht]
     \centering
    \begin{subfigure}[t]{0.32\textwidth}
        \raisebox{-\height}{\includegraphics[width=\textwidth]{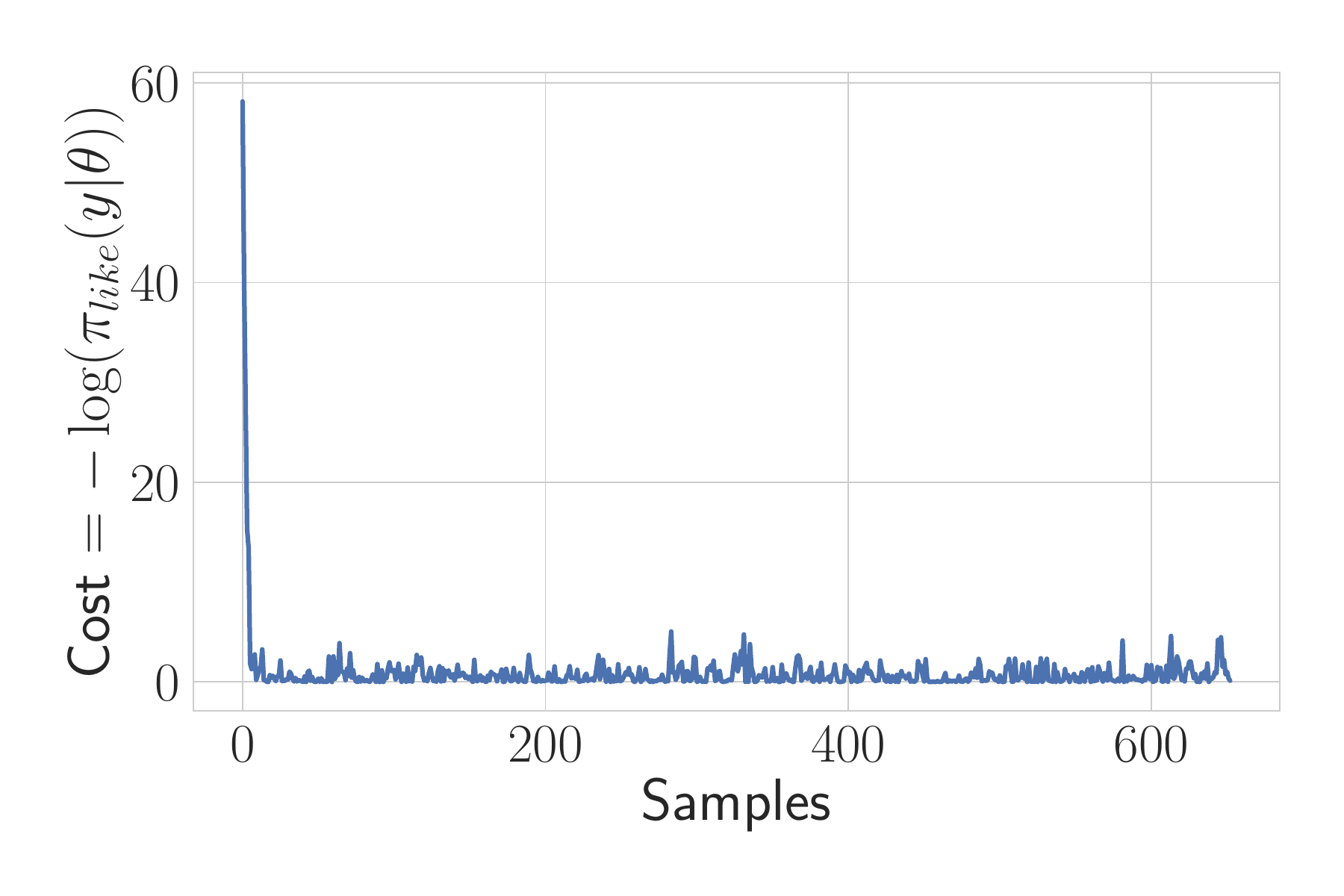}}
		\vspace{-5pt}
        \caption{Cost associated to the accepted samples.}
    \end{subfigure}
    \hfill
    \begin{subfigure}[t]{0.32\textwidth}
        \raisebox{-\height}{\includegraphics[width=\textwidth]{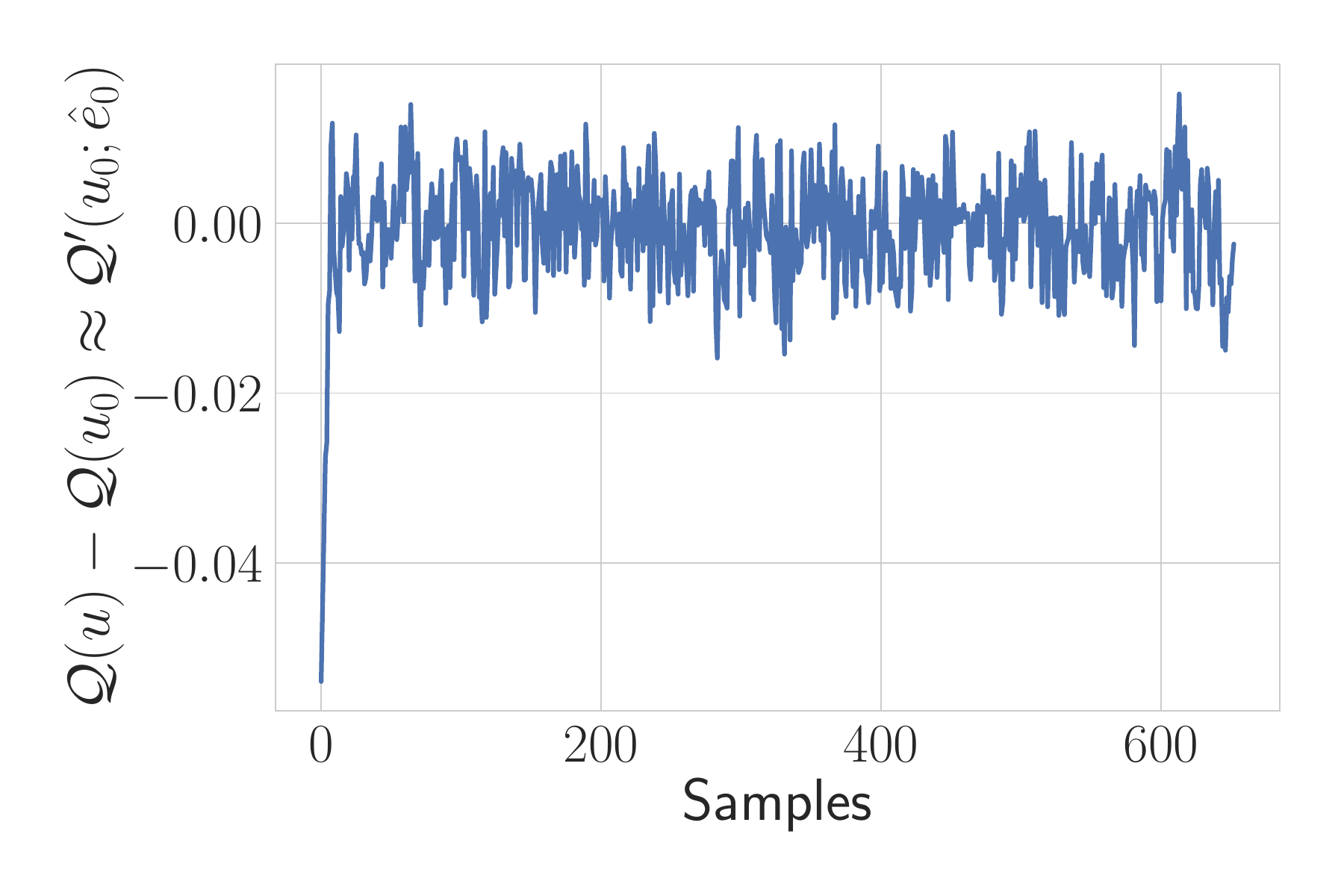}}
		\vspace{-5pt}
        \caption{QoI error $\calQ(u) - \calQ(u_0)$ using the estimate $\calQ'(u_0; \hat{e}_0)$.}
    \end{subfigure}
	\hfill
    \begin{subfigure}[t]{0.32\textwidth}
        \raisebox{-\height}{\includegraphics[width=\textwidth]{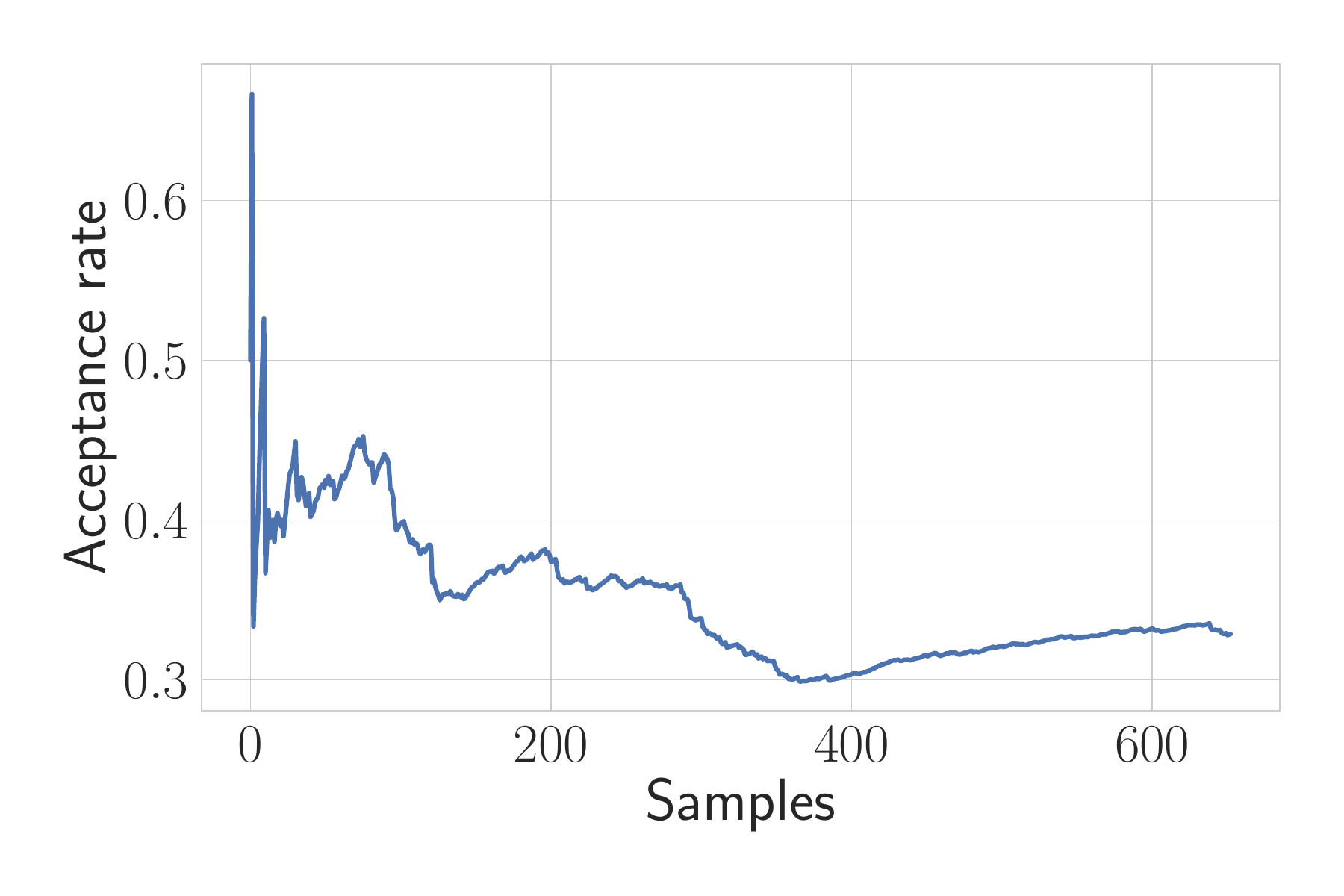}}
		\vspace{-5pt}
        \caption{Sample acceptance rate.}
    \end{subfigure}
    \caption{Results from one MCMC chain. In (A), the cost associated with the accepted samples is shown. In (B), the QoI error associated with the accepted samples is shown. The rate of sample acceptance is shown in (C). The cost and the QoI errors are seen to stabilize as more posterior samples are drawn indicating the convergence of the MCMC iterations. }\label{fig:mcmcStats2}
\end{figure}

\begin{figure}[ht]
	\centering
	\includegraphics[width=\textwidth]{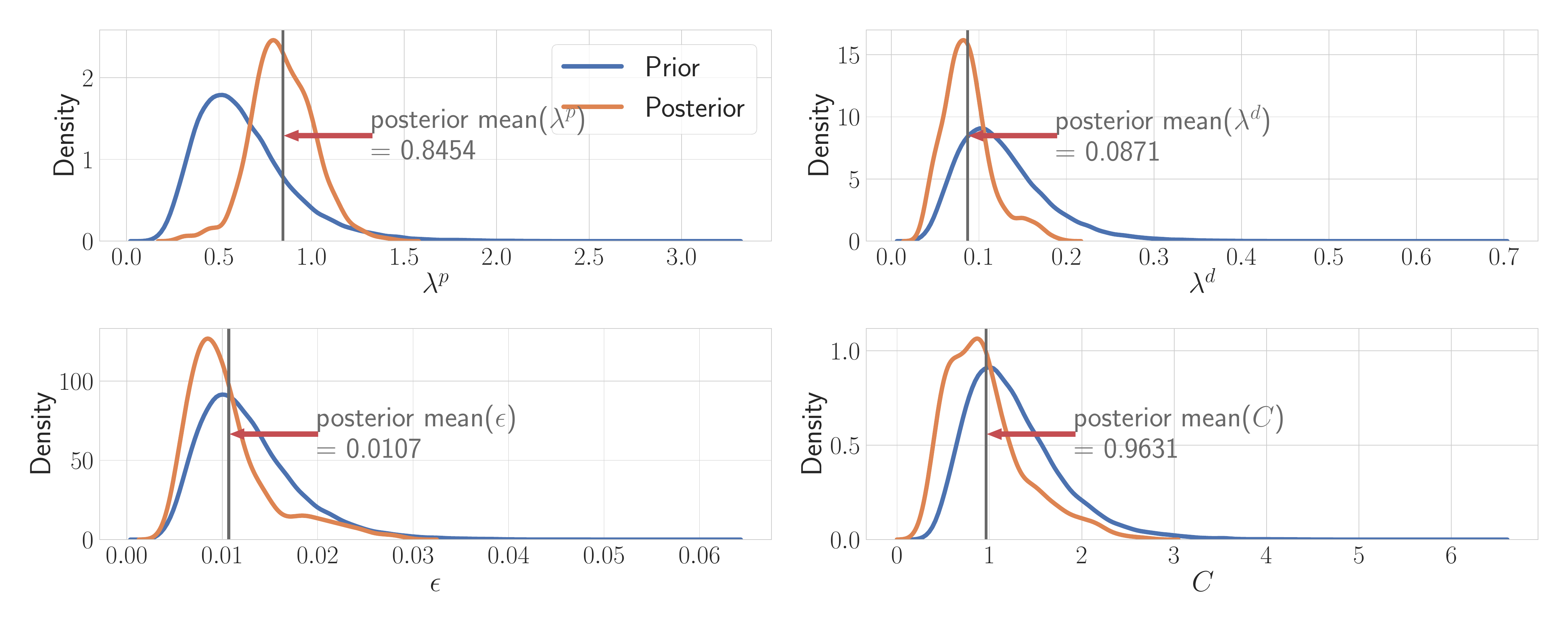}
	\vspace{-20pt}
	\caption{The prior and posterior probability density functions and mean values of posterior samples for parameters in the high-fidelity model.}\label{fig:comparePriorPost2}
\end{figure}

In \cref{fig:mcmcStats2}, the cost, the QoI error using the approximation $\calQ'(u_0; \hat{e}_0)$, and the sample acceptance rate for one of the chains are shown. 
In \cref{fig:comparePriorPost2}, the prior and posterior densities are shown. The mean and the standard deviation of the posterior samples are $\mu_{\btheta}  = (0.845, 0.087, 0.011, 0.963)$ and $\sigma_{\btheta} = (0.168, 0.028, 0.005, 0.435)$, respectively. The QoI error, $\calQ(u) - \calQ(u_0) \approx \calQ(\hat{e}_0)$, at $\btheta = \mu_{\btheta}$, is $-0.017$, i.e., $-1.4$ percent of $\calQ(u_0)$. 

%%%%%%%%%%%%%%%%%%%%%%%%%%%%%%%%%%%%%%%%%%%%%%%%%%%%%%%%%%%%%%%%%
%%%%%%%%%%%%%%%%%%%%%%%%%%%%%%%%%%%%%%%%%%%%%%%%%%%%%%%%%%%%%%%%%
\section{Conclusion}\label{s:conclusion}
In this work, a Bayesian approach for the calibration of parameters of a higher-fidelity PDE-based model with given parameters of a lower-fidelity model is presented. 
The technique works in the other direction as well; i.e., one can calibrate the lower-fidelity model if the higher-fidelity model is known. The central component of the proposed technique is the utilization of the goal-oriented error estimates. 
These estimates allow one to solve simpler linear problems for error components and to compute the error in quantities of interest approximately, but often with very good accuracy. The efficacy of the proposed approach is demonstrated by applying it to the two nontrivial examples.

While in this work the lower-fidelity model is assumed to be deterministic, there are useful scenarios where the lower-fidelity model is known but with some uncertainties in its parameters. 
The extension of the proposed approach to this case is straight-forward. 
Another interesting avenue for future work is the application of this approach to multifidelity Monte Carlo methods in which the models of varying fidelities are used to perform faster MCMC sampling.

Ongoing work includes the utilization of goal-oriented a-posteriori estimation in modeling the correction term in Bayesian inference problems in which the forward model is replaced with the surrogate model; it is possible to use the goal-oriented estimates to compute the error (error due to replacement of a model with the surrogate) cheaply, and use this error as a correction in a Bayesian inference framework. This work is in the same spirit as \cite{manzoni2016accurate}; however, rather than building the correction term offline as in \cite{manzoni2016accurate}, the correction can be computed on the fly as new parameters are sampled.

\section{Acknowledgments}
This work was supported by the U.S. Department of Energy, Office of Science, USA, Office of Advanced Scientific Computing Research, Mathematical Multifaceted Integrated Capability Centers (MMICCS), under Award Number DE-SC0019303.

%-------------------------------------------------------------------------%
%\bibliography{main}
%\bibliographystyle{siam}

%-------------------------------------------------------------------------%

\appendix
\section{Proof of Lemma 1}\label{s:proof}
\begin{proof}
The following Taylor expansion of a possibly nonlinear functional $\calQ: \calV \to \bbR$ with appropriate regularity is first noted. Let $u, u_0 \in \calV$ and $e_0 = u - u_0$, then
\begin{align}
\calQ(u) - \calQ(u_0) &= \int_0^1 \calQ'(u_0 + s e_0; e_0) \dd s \label{eq:FormTaylor1} \\
&= \calQ'(u_0; e_0) + \int_0^1 \calQ''(u_0 + s e_0; e_0, e_0) (1-s) \dd s \label{eq:FormTaylor2} \\
&= \frac{1}{2} \left[ \calQ'(u_0; e_0) + \calQ'(u; e_0) + \int_0^1 \calQ'''(u_0 + s e_0; e_0, e_0, e_0) (s-1)s \dd s \right] . \label{eq:FormTaylor3}
\end{align}
Suppose now that $(u_0, p_0)$ is any approximation of the solution $(u, p)$ of  problem \eqref{eq:fwdAdj}. Since $u$ satisfies $\calF(q) - \calB(u; q) = 0$, for any $q\in \calV$, the following holds,
\begin{align}\label{eq:Rest}
\calR(u_0; q) = \calF(q) - \calB(u_0; q) = \calF(q) - \calB(u; q) + \calB(u; q) - \calB(u_0; q) =  \calB(u; q) - \calB(u_0; q) .
\end{align}
Now, using \eqref{eq:FormTaylor2}, it is easy to show that
\begin{equation*}
 \calB(u; q) - \calB(u_0; q) = \calB'(u_0; e_0, q) + \int_0^1 \calB''(u_0 + s e_0; e_0, e_0, q) (1 - s) \dd s
\end{equation*}
and, therefore, from \eqref{eq:Rest}, it can be shown that
\begin{equation}\label{eq:RB1}
 \calB'(u_0; e_0, q) = \calR(u_0; q) - \int_0^1 \calB''(u_0 + s e_0; e_0, e_0, q) (1 - s) \dd s .
\end{equation}
Further, using \eqref{eq:FormTaylor3}, the following relation can also be obtained,
\begin{equation}\label{eq:Bdiff}
\begin{split}
 &\calB(u; q) - \calB(u_0; q) \\
&= \frac{1}{2} \left[ \calB'(u_0; e_0, q) + \calB'(u; e_0, q) + \int_0^1 \calB'''(u_0 + s e_0; e_0, e_0, e_0, q) (s-1)s \dd s \right] .
\end{split}
\end{equation}
Term $\calB'(u; e_0, q)$ can be estimated further, using \eqref{eq:FormTaylor2}, as follows
\begin{equation}
\begin{split}
\calB'(u; e_0, q) &= \calB'(u_0; e_0, q) + \calB'(u; e_0, q) - \calB'(u_0; e_0, q) \\
&= \calB'(u_0; e_0, q) + \calB''(u_0; e_0, e_0, q) + \int_0^1 \calB'''(u_0 + s e_0; e_0, e_0, e_0, q) (1 - s) \dd s .
\end{split}
\end{equation}
Combining above equation and \eqref{eq:Bdiff} with \eqref{eq:Rest} gives
\begin{align*}
\calR(u_0; q) =  \frac{1}{2} \left[ 2\calB'(u_0; e_0, q) + \calB''(u_0; e_0, e_0, q) + \int_0^1 \calB'''(u_0 + s e_0; e_0, e_0, e_0, q) ((1 - s) + (s-1)s) \dd s \right] .
\end{align*}
Thus, for any $q\in \calV$, it is shown that
\begin{equation}\label{eq:RB2}
\calB'(u_0; e_0, q) = \calR(u_0; q) -  \frac{1}{2} \calB''(u_0; e_0, e_0, q) -  \frac{1}{2}\int_0^1 \calB'''(u_0 + s e_0; e_0, e_0, e_0, q) (1 - s)^2 \dd s .
\end{equation}

Now consider $\bar{\calR}$. For any $v\in \calV$, 
\begin{align}\label{eq:Rbarest}
\bar{\calR}(u_0; v, p_0) &= \calQ'(u_0; v) - \calB'(u_0; v, p_0) \notag \\
&= \calQ'(u_0; v) - \calQ'(u; v) + \calQ'(u; v) - \calB'(u_0; v, p_0) \notag \\
&= \calQ'(u_0; v) - \calQ'(u; v)  + \calB'(u; v, p) - \calB'(u_0; v, p_0) \notag \\
&= - \left[ \calQ'(u; v) - \calQ'(u_0; v)\right] + \left[ \calB'(u; v, p) - \calB'(u_0; v, p) \right] + \left[\calB'(u_0; v, p) - \calB'(u_0; v, p_0) \right] \notag \\
&= - \left[ \calQ'(u; v) - \calQ'(u_0; v)\right] + \left[ \calB'(u; v, p) - \calB'(u_0; v, p) \right] + \calB'(u_0; v, \eps_0) ,
\end{align}
where, in the second step, $\calQ'(u; v) = \calB'(u; v, p)$ for all $v\in \calV$ is used.

It remains to establish \eqref{eq:AdjResFwdRes1}. From \eqref{eq:FormTaylor1}, one has
\begin{align*}
\calQ'(u; v) - \calQ'(u_0; v) &= \int_0^1 \calQ''(u_0 + s e_0; e_0, v) \dd s, \\
 \calB'(u; v, p) - \calB'(u_0; v, p) &= \int_0^1 \calB''(u_0 + s e_0; e_0, v, p) \dd s .
\end{align*}
Combining above with \eqref{eq:Rbarest} gives
\begin{align}\label{eq:RbarB1}
\bar{\calR}(u_0; v, p_0) = \calB'(u_0; v, \eps_0)  + \int_0^1 \left\{ \calB''(u_0 + s e_0; e_0, v, p)  - \calQ''(u_0 + s e_0; e_0, v) \right\} \dd s .
\end{align}
Taking $q = \eps_0$ in \eqref{eq:RB1} produces
\begin{equation*}
\calB'(u_0; e_0, \eps_0) = \calR(u_0; \eps_0) - \int_0^1 \calB''(u_0 + s e_0; e_0, e_0, \eps_0) (1 - s) \dd s .
\end{equation*}
Now, substituting $v=e_0$ in \eqref{eq:RbarB1}, taking into account the above relation, and noting $p - (1-s) \eps_0 = p_0 + s \eps_0$ gives \eqref{eq:AdjResFwdRes1}.

The relation \eqref{eq:AdjResFwdRes2} can be derived from \eqref{eq:Rbarest} and previous estimates. Applying \eqref{eq:FormTaylor2} to $\calQ'$ and $\calB'$ gives
\begin{align*}
\calQ'(u; v) - \calQ'(u_0; v) &= \calQ''(u_0; e_0, v) + \int_0^1 \calQ'''(u_0 + s e_0; e_0, e_0, v) (1 - s) \dd s, \\
 \calB'(u; v, p) - \calB'(u_0; v, p) &= \calB''(u_0; e_0, v, p) + \int_0^1 \calB'''(u_0 + s e_0; e_0, e_0, v, p) (1 - s) \dd s .
\end{align*}
Substituting the above into \eqref{eq:Rbarest} gives
\begin{align}\label{eq:RbarB2}
\bar{\calR}(u_0; v, p_0) &= \calB'(u_0; v, \eps_0) - \calQ''(u_0; e_0, v) + \calB''(u_0; e_0, v, p) \notag \\
&\quad + \int_0^1 \left\{\calB'''(u_0 + s e_0; e_0, e_0, v, p) - \calQ'''(u_0 + s e_0; e_0, e_0, v) \right\} (1 - s) \dd s 
\end{align}
Taking $q = \eps_0$ in \eqref{eq:RB2} produces
\begin{equation*}
\calB'(u_0; e_0, \eps_0) = \calR(u_0; \eps_0) -  \frac{1}{2} \calB''(u_0; e_0, e_0, \eps_0) -  \frac{1}{2}\int_0^1 \calB'''(u_0 + s e_0; e_0, e_0, e_0, \eps_0) (1 - s)^2 \dd s .
\end{equation*}
Substituting $v=e_0$ in \eqref{eq:RbarB2} and using above relation, it can be shown that
\begin{align*}
&\bar{\calR}(u_0; e_0, p_0) \notag \\
&= \calR(u_0; \eps_0) -  \frac{1}{2} \calB''(u_0; e_0, e_0, \eps_0) - \calQ''(u_0; e_0, e_0) + \calB''(u_0; e_0, e_0, p)  + \int_0^1 \left\{\calB'''(u_0 + s e_0; e_0, e_0, e_0, p) \right. \notag \\
&\qquad \left. - \calQ'''(u_0 + s e_0; e_0, e_0, e_0) -\frac{1}{2} \calB'''(u_0 + s e_0; e_0, e_0, e_0, \eps_0) (1 - s) \right\} (1 - s) \dd s  \notag \\
&=  \calR(u_0; \eps_0)- \calQ''(u_0; e_0, e_0) + \calB''(u_0; e_0, e_0, p_0)  +  \frac{1}{2} \calB''(u_0; e_0, e_0, \eps_0) \notag \\
&\quad \int_0^1 \left\{\calB'''\left(u_0 + s e_0; e_0, e_0, e_0, p - \frac{1-s}{2}\eps_0\right) - \calQ'''(u_0 + s e_0; e_0, e_0, e_0)\right\} (1 - s) \dd s .
\end{align*}
This completes the poof of the lemma.
\end{proof}

\section{Numerical discretization of nonlinear and linear transient problems}\label{s:numdisc}
Suppose $V_h \subset H^1(\Omega)$ is the finite element function space spanned by bi-linear interpolation functions on a quadrilateral mesh $\calT$ of domain $\Omega$. To solve the nonlinear problem numerically, a first-order semi-implicit time marching scheme is employed. Picard iteration is used to solve the nonlinear problem at each time step. Let $u_n\in V_h$ be the solution at time $t_n = n \Delta t$, $\Delta t$ being the time step, and let $u^k$, for $k = 0, 1, 2, ..,$, being the \kth{k} iterative solution at time $t_{n+1}$. Given $u_n$ and \kth{k} iterative solution $u^k$, the variational problem to compute the next iterative solution $u^{k+1}$ at time $t_{n+1}$ is written as
\begin{equation}\label{eq:iterscheme}
\begin{split}
&\weakDot{\frac{u^{k+1} - u_n}{\Delta t}}{v} + \weakDot{\epsilon \nabla u^{k+1}}{\nabla v} + \weakDot{\lambda^d u^{k+1}}{v} \\
&+ \weakDot{\Psi'_{im}(u^{k+1}) + \Psi'_{ex}(u^{k})}{v} - \weakDot{\lambda^p u^{k+1} (1 - u^k) f}{v} = 0, \qquad\qquad \forall v \in V_h,
\end{split}
\end{equation}
where $\Psi$ has been decomposed in two parts, $\Psi_{im} = \Psi_{im}(r)$ and $\Psi_{ex} = \Psi_{ex}(r)$; $\Psi_{im}$ is quadratic in $r$ (so $\Psi'_{im}$ is linear) and  therefore it is treated implicitly in the above scheme. On the other hand, $\Psi_{ex}$ is a convex function and is treated explicitly. These are defined as
\begin{equation}
\Psi_{im}(r) = \frac{3C}{2} r^2, \qquad \Psi_{ex}(r) = \frac{C}{2} (2r^4 - 4r^3 -  r^2). 
\end{equation}
Having solved for $u^{k+1}$ using \eqref{eq:iterscheme}, the next iterative solution $u^{k+2}$ is only sought if the error in the current and old iterative solutions is above a specified tolerance. I.e., if the following,
\begin{equation}
||u^{k+1} - u^{k}||_{L^2(\Omega)} < \mathrm{tol} := 10^{-10}, 
\end{equation}
holds, the iteration is stopped and the solution $u_{n+1}$ at time $t_{n+1}$ is taken as $u_{n+1} = u^{k+1}$. 

To solve the linear problem for an approximate error $\hat{e}_0$, see \eqref{eq:errWeakForm1}$_1$, a first-order implicit time marching scheme is used. Given a field $u_0 = u_0(t, \bx)$ and the solution $\hat{e}_{0_n}$ at time $t_n$, the variational problem to compute $\hat{e}_{0_{n+1}}$ at time $t_{n+1}$ is written as
\begin{equation}
\begin{split}
&\weakDot{\frac{\hat{e}_{0_{n+1}} - \hat{e}_{0_n}}{\Delta t}}{v} + \weakDot{\epsilon \nabla \hat{e}_{0_{n+1}} }{\nabla v} + \weakDot{\lambda^d \hat{e}_{0_{n+1}} }{v} \\
&\quad + \weakDot{\Psi''(u_{0_{n+1}}) \hat{e}_{0_{n+1}}}{v} - \weakDot{\lambda^p(1 - 2 u_{0_{n+1}}) \hat{e}_{0_{n+1}} f}{v} \\
&= -\weakDot{\frac{u_{0_{n+1}} - u_{0_n}}{\Delta t}}{v} - \weakDot{\epsilon \nabla u_{0_{n+1}} }{\nabla v} - \weakDot{\lambda^d u_{0_{n+1}} }{v} \\
&\quad -\weakDot{\Psi'(u_{0_{n+1}})}{v} + \weakDot{\lambda^p u_{0_{n+1}} (1 - u_{0_{n+1}}) f }{v}, 
\end{split}
\end{equation}
for all $v \in V_h$.

\end{document}

%% file: preamble.tex
\usepackage{mathtools,bm,comment,amssymb,enumitem,cite,xcolor,array,float,subcaption,float,pgfplots}
\usepackage[capitalise, nameinlink]{cleveref}
\usepackage[english]{babel}
\usepackage[margin=1in]{geometry}
\usepackage{amsaddr}
\usepackage[linesnumbered,ruled]{algorithm2e}
\usepackage[noend]{algpseudocode}
\usepackage[locale=DE]{siunitx}
\usepackage{textcomp}
\usetikzlibrary{shapes,arrows}

\usepackage{url}

\usepackage{graphicx}

\usepackage{tcolorbox}

\usepackage{amsthm}
\theoremstyle{remark}

\numberwithin{equation}{section}
\crefname{subsection}{Subsection}{Subsections}
\crefformat{equation}{(#2#1#3)}
\crefformat{enumi}{(#2#1#3)}
\crefname{figure}{Figure}{Figures}

\newcommand{\eps}{\varepsilon}

\newcommand{\dd}{\,\textup{d}}

\renewcommand{\t}{\tilde}
\newcommand{\p}{\partial}

\newtheorem{remark}{Remark}

\newtheorem{lemma}{Lemma}
\newtheorem{theorem}{Theorem}

\DeclareMathOperator{\sign}{sign}

\newcommand{\bbR}{\mathbb{R}}

\newcommand{\normX}[2]{{\left\lVert #1 \right\rVert}_{#2}}

\newcommand{\leftcr}{\left\{}
\newcommand{\rightcr}{\right\}}

\newcommand{\calF}{{\mathcal{F}}}

\newcommand{\calA}{{\mathcal{A}}}
\newcommand{\calB}{{\mathcal{B}}}
\newcommand{\calN}{{\mathcal{N}}}
\newcommand{\calT}{{\mathcal{T}}}
\newcommand{\calU}{{\mathcal{U}}}
\newcommand{\calV}{{\mathcal{V}}}

\newcommand{\calQ}{{\mathcal{Q}}}

\newcommand{\calR}{{\mathcal{R}}}

\newcommand{\bn}{\bm{n}}

\newcommand{\bx}{\bm{x}}

\newcommand{\btheta}{\bm{\theta}}

\newcommand{\weakDot}[2]{\left({#1}, {#2}\right)}

\newcommand{\dualDot}[2]{\left< {#1}, {#2}\right>}

\newcommand{\kth}[1]{${#1}^{\text{th}}$}

\newcommand{\sqbrack}[1]{\left[ #1 \right]}
\newcommand{\cubrack}[1]{\left\{ #1 \right\}}
\newcommand{\orbrack}[1]{\left( #1 \right)}

\newenvironment{tbox}{\begin{tcolorbox}[colback=white]}{\end{tcolorbox}}

\newcommand{\formBExpand}[2]{\int_\Omega \left\{ \kappa (1+ #1^2) \nabla #1 \cdot \nabla #2  + \alpha #1 #2\right\}\dd \bx}

\newcommand{\formBDerExpand}[3]{\int_\Omega \left\{\kappa (1+#1^2) \nabla #2 \cdot \nabla #3 + 2\kappa #1 #2  \nabla #1\cdot \nabla #3 + \alpha #2 #3\right\} \dd \bx}

\usepackage{enumitem}
\setlist[itemize]{left=5pt}
\setlist[enumerate]{left=5pt}
\setlist[description]{left=5pt}

\usepackage[ruled]{algorithm2e}
\usepackage{algcompatible}
\algblockdefx[Hloop]{Hloop}{EndHloop}[1][]{\textbf{hpx::parallel::for$\_$loop} #1}{\textbf{end parallel for}}
\algblockdefx[Func]{Func}{EndFunc}[1][]{\textbf{Function }}{\textbf{end Function}}

\usepackage{listings}

\usepackage{xcolor}
\definecolor{mygray}{rgb}{0.1, 0.1, 0.8}

\makeatletter
\@namedef{subjclassname@2020}{\textup{2020} Mathematics Subject Classification}
\makeatother

\makeatletter
\renewcommand{\email}[2][]{%
	\ifx\emails\@empty\relax\else{\g@addto@macro\emails{,\space}}\fi%
	\@ifnotempty{#1}{\g@addto@macro\emails{\textrm{(#1)}\space}}%
	\g@addto@macro\emails{#2}%
}
\makeatother

\makeatletter
\pgfmathdeclarefunction{erf}{1}{%
	\begingroup
	\pgfmathparse{#1 > 0 ? 1 : -1}%
	\edef\sign{\pgfmathresult}%
	\pgfmathparse{abs(#1)}%
	\edef\x{\pgfmathresult}%
	\pgfmathparse{1/(1+0.3275911*\x)}%
	\edef\t{\pgfmathresult}%
	\pgfmathparse{%
		1 - (((((1.061405429*\t -1.453152027)*\t) + 1.421413741)*\t 
		-0.284496736)*\t + 0.254829592)*\t*exp(-(\x*\x))}%
	\edef\y{\pgfmathresult}%
	\pgfmathparse{(\sign)*\y}%
	\pgfmath@smuggleone\pgfmathresult%
	\endgroup
}
\makeatother

\let\originalleft\left
\let\originalright\right
\renewcommand{\left}{\mathopen{}\mathclose\bgroup\originalleft}
\renewcommand{\right}{\aftergroup\egroup\originalright}

\makeatletter

\renewcommand\subsubsection{\@secnumfont}{\bfseries}%
\renewcommand\subsubsection{\@startsection{subsubsection}{3}
  \z@{.5\linespacing\@plus.7\linespacing}{-.5em}%
  {\normalfont\bfseries}}
  
\renewcommand\paragraph{\@startsection{paragraph}{4}%
  \z@\z@{-\fontdimen2\font}%
  {\normalfont\bfseries}}

\makeatother